\documentclass[a4paper]{article}

\usepackage[T1]{fontenc}
\usepackage[a4paper, total={6in, 8.7in}]{geometry}
\usepackage{amssymb,amsmath}
\usepackage{dsfont}
\usepackage{algorithm}
\usepackage{algpseudocode}
\usepackage{amsthm}

\DeclareMathOperator*{\argmin}{arg\,min}
\DeclareMathOperator*{\argmax}{arg\,max}
\usepackage{amsmath,amssymb,amsthm,amsfonts}
\usepackage{nicefrac}

\begin{document}

\newtheorem{thm}{Theorem}[section]
\newtheorem{asm}{Assumption}[section]
\newtheorem{lemma}[thm]{Lemma}
\newtheorem{cor}[thm]{Corollary}
\newtheorem{dfn}{Definition}[section]
\newtheorem{exm}{Example}[section]
\newtheorem{rmk}{Remark}[section]
\allowdisplaybreaks

\title{\textbf{Consistency of Variational Bayes Inference for Estimation and Model Selection in Mixtures}}
\author{Badr-Eddine Ch\'erief-Abdellatif \& Pierre Alquier\footnote{Both authors gratefully acknowledge financial support from GENES and  by the French National Research Agency (ANR) under the grant Labex Ecodec (ANR-11-LABEX-0047). The second author gratefully acknowledges financial support from  the  research  programme New Challenges for New Data from LCL and GENES, hosted by the Fondation du Risque.}
\\
\small{CREST, ENSAE, Universit\'e Paris Saclay}}
\date{\today}

\maketitle

\begin{abstract}
\noindent Mixture models are widely used in Bayesian statistics and machine learning, in particular in computational biology, natural language processing and many other fields. Variational inference, a technique for approximating intractable posteriors thanks to optimization algorithms, is extremely popular in practice when dealing with complex models such as mixtures. The contribution of this paper is two-fold. First, we study the concentration of variational approximations of posteriors, which is still an open problem for general mixtures, and we derive consistency and rates of convergence. We also tackle the problem of model selection for the number of components: we study the approach already used in practice, which consists in maximizing a numerical criterion (the Evidence Lower Bound). We prove that this strategy indeed leads to strong oracle inequalities. We illustrate our theoretical results by applications to Gaussian and multinomial mixtures.
\end{abstract}

\section{Introduction}

\vspace{0.2cm}

This paper studies the statistical properties of variational inference as a tool to tackle two problems of interest: estimation and model selection in mixture models. Mixtures are often used for modelling population heterogeneity, leading to practical clustering methods \cite{Bouveyron2014,mcnicholas2016}. Moreover they have enough flexibility to approximate accurately almost every density \cite{ApproximatingDensitiesWithMixtures,kruijer2010adaptive}. Mixtures are used in many various areas such as computer vision \cite{MixturesInComputerVision}, genetics \cite{MixturesinGenetics}, economics \cite{MixturesInEconomics}, transport data analysis \cite{carel2017simultaneous} and others. We refer the reader to \cite{handbook} for an account of the recent advances on mixtures. The most famous procedure for mixture density estimation in the frequentist literature is probably Expectation-Maximization \cite{EM}, a maximum-likelihood algorithm that yields increasingly higher likelihood. At the same time, the Bayesian paradigm has raised great interest among researchers and practitioners, especially through the Variational Bayes (VB) framework which aims at maximizing a quantity referred to as Evidence Lower Bound on the marginal likelihood (ELBO). Variational Bayes inference is a useful tool for approximating  intractable posteriors. It is known to work well in practice for mixture models: one of the most recent survey on VB \cite{Review} chooses mixtures as an example of choice to illustrate the power of the method. Moreover  \cite{Review} states: "the [evidence lower] bound is a good approximation of the marginal likelihood, which provides a basis for selecting a model. Though this sometimes works in practice, selecting based on a bound is not justified in theory". The main contribution of this paper is to prove that VB is consistent for estimation in mixture models, and that the ELBO maximization strategy used in practice is consistent for model selection. Thus we solve the question raised by \cite{Review}.
\vspace{0.2cm}

Variational Bayes is a method for computing intractable posteriors in Bayesian statistics and machine learning. Markov Chain Monte Carlo (MCMC) algorithms remain the most widely used methods in computational Bayesian statistics. Nevertheless, they are often too slow for practical uses when the dataset is very large. A more and more popular alternative consists in finding a deterministic approximation of the target distribution called Variational Bayes approximation. The idea is to minimize the Kullback-Leibler divergence of a tractable distribution $\rho$ with respect to the posterior, which is also equivalent to maximizing the ELBO. This optimization procedure is much faster and efficient than MCMC sampling with numerous applications in different fields: matrix completion for collaborative filtering \cite{MatrixCompletionAlquierCottet}, computer vision \cite{Vision}, computational biology \cite{Biology} and natural language processing \cite{hoffman2013stochastic}, to name a few prominent examples.

\vspace{0.2cm}

However, variational inference is mainly used for its practical efficiency and only little attention has been put in the literature towards theoretical properties of the VB approximation until very recently. In \cite{alquier2016properties} the properties of variational approximations of Gibbs distributions used in machine learning are derived. The results are essentially valid for bounded loss functions, which makes them difficult to use beyond the problem of supervised classification. Based on some technical advances from \cite{bhattacharya2016bayesian}, \cite{Tempered} removed the boundedness assumption in \cite{alquier2016properties}, allowing to study more general statistical models. In \cite{Plage}, the authors extended the range of models covered by \cite{alquier2016properties}. This allowed them to study mixture of Gaussian distributions as an example. Many questions are still left unanswered: model selection, and the estimation of mixture of non-Gaussian distributions. For example mixture of multinomials are widely used in practice \cite{carel2017simultaneous}, as well as more intricated examples such as nonparametric mixtures \cite{gassiat2018efficient}. Note that all the results in \cite{bhattacharya2016bayesian,Tempered,Plage} are limited to so-called tempered posteriors, that is, where the likelihood is taken to some power $\alpha$. Still, the use of tempered posteriors is highly recommended by many authors as a way to overcome model misspecification, see \cite{grunwaldmisspecifiation} and the references therein. Indeed some results in \cite{Tempered} are valid in a misspecified setting. Alternative approaches were developed to study VB: \cite{BVM} established Bernstein-von-Mises type theorems on the variational approximation of the posterior. They provide very interesting results for parametric models but it is unclear whether these results can be extended to model selection or misspecified case. More recently, \cite{Chicago} succeeded in adapting the now classical results of \cite{ghosal2000convergence} to Variational Bayes and showed that a slight modification in the three classical "prior mass and testing conditions" leads to the convergence of their variational approximations, again under the assumption that the model is true. With respect to these works, our contribution is a complete study of the consistency of VB for mixtures of general distributions. In particular, we explicit independent conditions on the prior on the weights, and on the prior on the parameters of the components. The study is done in the case $\alpha<1$ which allows to prove results in the misspecified case.

\vspace{0.2cm}

The other point addressed in this paper is model selection. This is a natural question which can be interpreted in this context as the determination of the number of components of the mixture. This point is crucial: indeed, too many components can lead to estimates with too large variances whereas with too few components, we may obtain mixtures which are not able to fit the data properly. This is a common issue and a lot of statisticians worked on this question. In the literature, criteria such as AIC \cite{AIC} and BIC \cite{BIC} are popular. It is well known that in some collections of models, AIC optimizes the prediction ability while BIC recovers with high probability the true model (when there is one). These two objectives are not compatible in general \cite{yang2005can}. Anyway, these results depend on asssumptions that are not satisfied by mixtures. It seems thus more natural to develop criteria suited to a given objective. For example, \cite{biernacki1999improvement} proposed a procedure to select a number of components that is the most relevant for clustering. A non-asymptotic theory of penalization has been developed during the last two decades using oracle inequalities \cite{Massart}. In the wake of those works, our paper studies mixture model selection based on the ELBO criterion. We prove a general oracle inequality. This result establishes the consistency of ELBO maximization when the primary objective is the estimation of the distribution of the data.

\vspace{0.2cm}

The rest of this paper is organized as follows. In Section \ref{sec:notations} we introduce the background and the notations that will be adopted. Consistency of the Variational Bayes for estimation in a mixture model is studied in Section \ref{sec:main-res}. First, we give the general results under a "prior mass" assumption, as well as a general form for the algorithm to compute the VB approximation (Subsection \ref{subsec:pacbayes}). We then apply these results to mixtures of multinomials (Subsection \ref{sec:main-res:2}) and Gaussian mixtures (Subsection \ref{sec:main-res:3}). In each case, we provide a rate of convergence of VB and discuss its numerical implementation. We extend the setting to the misspecified case in Subsection \ref{sec:main-res:4}. Finally, we address the issue of selecting based on the ELBO in Section \ref{sec:selection}. We discuss possible extensions in Section~\ref{section:conclusion}, while Section \ref{sec:proofs} is dedicated to the proofs.

\vspace{0.2cm}

\section{Background and notations}
\label{sec:notations}

\vspace{0.2cm}

Let us introduce the notations and the framework we adopt in this paper. We observe in a measurable space $\big( \mathbb{X},\mathcal{X} \big)$ a collection of $n$ i.i.d. random variables $X_1$,...,$X_n$ sampled from a probability distribution $P^0$. We put $(X_1,...,X_n) = X_1^n$. The goal is to estimate the generating distribution $P^0$ of the $X_i$'s by a $K$-components mixture model. We will study the (frequentist) properties of variational approximations of the posterior. The extension to selection of the number of components is also tackled in this paper, but we will first deal with a fixed $K$. We introduce a collection of distributions $\{ Q_{\theta}/ \theta \in \Theta \}$ indexed by a parameter space $\Theta$ from which we will take the different components of our mixture model. We assume that for each $\theta \in \Theta$, the probability distribution $Q_\theta$ is dominated by a reference measure $\mu$ and that the density $q_\theta=\frac{dQ_\theta}{d\mu}$ is such that the map $(x,\theta) \rightarrow q_\theta(x)$ is $\mathcal{X} \times \mathcal{T}$-measurable, $\mathcal{T}$ being some sigma-algebra on $\Theta$. Unless explicitly stated otherwise, all the distributions that will be considered in this paper will be characterized by their density with respect to the dominating measure $\mu$. We can now consider the statistical mixture model of $K \geq 1$ components defined as:
$$
\left\lbrace \sum_{j=1}^K p_j Q_{\theta_j} \hspace{0.2cm} \big/ \hspace{0.2cm} \theta_j \in \Theta \hspace{0.2cm} \textnormal{for} \hspace{0.2cm} j=1,...,K , \hspace{0.2cm} p=(p_1,...,p_K) \in \mathcal{S}_K \right\rbrace
$$ 
where $\mathcal{S}_K=\big\{ p=(p_1,...,p_K) \in \mathbb{R}^K \big/ \hspace{0.2cm} p_j \geq 0  \hspace{0.2cm} \textnormal{for} \hspace{0.2cm} j=1,...,K \hspace{0.2cm} \textnormal{and} \hspace{0.2cm} \sum_{j=1}^K p_j = 1 \big\}$ is the $K-1$ dimensional simplex. We will write $\theta=(p_1,...,p_K,\theta_1,...\theta_K) \in \Theta_K$ for short, where $p \in \mathcal{S}_K$, $\theta_j \in \Theta$ for $j=1,...,K$ and $\Theta_K=\mathcal{S}_K \times \Theta^K$. The mixture corresponding to parameter $\theta=(p_1,...,p_K,\theta_1,...,\theta_K)$ will be denoted $P_{\theta}:=\sum_{j=1}^K p_j Q_{\theta_j}$.

\vspace{0.2cm}
First, we consider the well-specified case, assuming that the true distribution belongs to the $K$-components mixture model. Thus, we define the true distribution $P^0$ from which data are sampled: 
$$ 
X_1,...,X_n \sim \sum_{j=1}^K p^0_j Q_{\theta^0_j} \hspace{0.2cm} \textnormal{with} \hspace{0.2cm} \theta^0_j \in \Theta \hspace{0.2cm} \textnormal{for} \hspace{0.2cm} j = 1,...,K \hspace{0.2cm} \textnormal{and} \hspace{0.2cm} p^0 \in \mathcal{S}_K.
$$
Hence, we want to estimate the true distribution $P_{\theta^0}$ using a Bayesian approach. Therefore, we define a prior $\pi=\pi_p \bigotimes_{j=1}^{K}\pi_j$ on $\theta$, $\pi_p \in \mathcal{M}_1^+(\mathcal{S}_K)$ being a probability distribution on some measurable space $(\mathcal{S}_K,\mathcal{A})$, and each $\pi_j \in \mathcal{M}_1^+(\Theta)$ a probability distribution on the measurable space $(\Theta,\mathcal{T})$. We will also consider in this paper the misspecified case where the true distribution does not belong to our statistical model i.e. is not necessarily a mixture, but the specific notations and framework will be described later. 

\vspace{0.2cm}
Let us introduce some notations. The likelihood will be denoted by $L_n$ and the log-likelihood by $\ell_n$, that is, for any $\theta=(p_1,...,p_K,\theta_1,...\theta_K)$,
\[
L_n(\theta) = \prod_{i=1}^n \sum_{j=1}^K p_j q_{\theta_j} \hspace{0.1cm}, \hspace{0.2cm}
\ell_n(\theta) = \sum_{i=1}^n \log\bigg({\sum_{j=1}^K p_j q_{\theta_j}}\bigg).
\]
The negative log-likelihood ratio $r_n$ between two distributions $P$ and $R$ is given by
\[
r_n(P,R) = \sum_{i=1}^n \log\bigg( \frac{R(X_i)}{P(X_i)} \bigg)
\]
(note that $r_n(\theta,\theta')$ is used by many authors instead of $r_n(P_{\theta},P_{\theta'})$ but our notation is more convenient for the extension to the misspecified case). The Kullback-Leibler (KL) divergence between two probability distributions $P$ and $R$ is given by
$$
\mathcal{K}(P,R) = \begin{cases}
\int \log\left( \frac{dP}{dR} \right) dP \hspace{0.2cm} \text{if $R$ dominates $P$}, \\
+ \infty \hspace{0.2cm} \text{otherwise.}
\end{cases}
$$
If some measure $\lambda$ dominates both $P$ and $R$ distributions represented here by their densities $f$ and $g$ with respect to this measure, we have
\[
\mathcal{K}(P,R) = \int f \log\left( \frac{f}{g} \right) d\lambda
\]
and we will use $K(P,R)$ or $\mathcal{K}(f,g)$ to denote this quantity, depending on the context.

\vspace{0.2cm}
 We also remind that the $\alpha$-Renyi divergence between $P$ and $R$,
$$
D_\alpha(P,R) = \begin{cases}
\frac{1}{\alpha-1} \log \int \left( \frac{dP}{dR} \right) ^{\alpha-1} dP \hspace{0.2cm} \text{if $R$ dominates $P$}, \\
+ \infty \hspace{0.2cm} \text{otherwise.}
\end{cases}
$$
When for some $\lambda$ we have $f = \frac{dP}{d\lambda}$ and $g=\frac{dR}{d\lambda}$,
\[
D_\alpha(P,R) = D_\alpha(f,g) = \frac{1}{\alpha-1} \log \int f^\alpha g^{1-\alpha} d\lambda.
\]

\vspace{0.2cm}
Some useful properties of Renyi divergences can be found in \cite{van2014renyi}. In particular, the Renyi divergence between two probability distributions $P$ and $R$ can be related to the classical total variation $TV$ and Hellinger $H$ distances respectively defined as $TV(P,R)=\frac{1}{2} \int |dP-dR|$ and $H(P,R)^2={ \frac{1}{2}\int (\sqrt{dP}-\sqrt{dR})^2 }={1-e^{-\frac{1}{2}D_{1/2}(P,R)}}$ through:
$$
TV(P,R)^2 \leq 2 H(P,R)^2 \leq D_{1/2}(P,R) \text{ and } D_{\alpha}(P,R) \xrightarrow[\alpha \rightarrow 1]{\nearrow} \mathcal{K}(P,R).
$$

\vspace{0.2cm}
The tempered Bayesian posterior $\pi_{n,\alpha}(.|X_1^n)$, which is our target here, is defined for $0<\alpha \leq 1$ by
\[
\pi_{n,\alpha}(d\theta|X_1^n) =\frac{e^{-\alpha r_n(P_\theta,P^0)}\pi(d\theta)}{\int e^{-\alpha r_n(P_\phi,P^0)}\pi(d\phi)} \propto L_{n}(\theta)^\alpha \pi(d\theta)
\]
(it is also referred to as fractional posterior, for example in \cite{bhattacharya2016bayesian}).
Note that when $\alpha=1$, then we recover the "true" Bayesian posterior, but the case $\alpha<1$ has many advantages: it is often more tractable from a computational perspective \cite{neal1996sampling,behrens2012tuning}, it is consistent under less stringent assumptions than required for $\alpha=1$ \cite{bhattacharya2016bayesian} and it is more robust to misspecification \cite{grunwaldmisspecifiation}.

\vspace{0.2cm}
We are now in position to define the VB approximation $\tilde{\pi}_{n,\alpha}(.|X_1^n)$ of the tempered posterior with respect to some set of distributions $\mathcal{F}$: it is the projection, with respect to the Kullback-Leibler divergence, of the tempered posterior onto the mean-field variational set $\mathcal{F}$,
\begin{equation*}
\tilde{\pi}_{n,\alpha}(.|X_1^n) = \argmin_{\rho \in \mathcal{F}} \mathcal{K}\bigg(\rho,\pi_{n,\alpha}(.|X_1^n)\bigg).
\end{equation*}
The mean-field approximation is very popular in the Variational Bayes literature. It is based on a decomposition of the space of parameters $\Theta_K$ as a product. Then $\mathcal{F}$ consists in compatible product distributions. Here, a natural choice \cite{Review} is $\Theta_K=\mathcal{S}_K \times \Theta \times \dots \times \Theta $ and
$$
\mathcal{F}=\bigg\{\rho_p \bigotimes_{j=1}^{K} \rho_j / \rho_p \in \mathcal{M}_1^+(\mathcal{S}_K), \hspace{0.1cm} \rho_j \in \mathcal{M}_1^+(\Theta) \hspace{0.1cm} \forall j=1,...,K\bigg\}.
$$
We will work on this particular set in the following and we will often use $\rho$ instead of $\rho_p \bigotimes_{j=1}^{K} \rho_j$ to ease notation.

\vspace{0.2cm}
We end this section by recalling Donsker and Varadhan's variational formula. Refer for example to \cite{MR2483528} for a proof (Lemma 1.1.3).

\begin{lemma}
\label{thm-dv}
For any probability $\lambda$ on some measurable space $(\textbf{E},\mathcal{E})$ and any measurable function
  $h: \textbf{E} \rightarrow \mathbb{R}$ such that $\int{\rm e}^h  \rm{d}\lambda < \infty$,
  \begin{equation*}
    \log\int {\rm e}^h \mathrm{d}\lambda = \underset{\rho \in \mathcal{M}_{1}^+(\textbf{E})}{\sup} \bigg\{ \int h \mathrm{d}\rho - \mathcal{K}(\rho,\lambda) \bigg\},
  \end{equation*}
  with the convention
  $\infty-\infty =-\infty$. Moreover, if $h$ is upper-bounded on the
  support of $\lambda$, then the supremum  on the right-hand
  side is reached by the distribution of the form:
  \begin{equation*}
    \lambda_h(d\beta) =
    \frac{{\rm e}^{h(\beta)} }{\int{\rm e}^h \mathrm{d}\lambda} \lambda(\mathrm{d}\beta).
  \end{equation*}
\end{lemma}
This technical lemma is one of the main ingredients for the proof of our results, but it is also very helpful to understand variational approximations. Indeed,
for $\textbf{E}=\Theta_K$ and using the definition of ${\pi}_{n,\alpha}(.|X_1^n)$, we get:
$$
 \pi_{n,\alpha}(\cdot|X_1^n)  = \underset{\rho \in \mathcal{M}^1_+(\Theta_K)}{\arg\min} \, \left\{
 \alpha \int r_n(P_{\theta},P^0) \rho({\rm d}\theta) + \mathcal{K}(\rho,\pi)
 \right\}
$$
and simple calculations give
\begin{align} \nonumber
  \tilde{\pi}_{n,\alpha}(\cdot|X_1^n) & = \underset{\rho \in \mathcal{F}}{\arg\min} \, \bigg\{
 \alpha \int r_n(P_{\theta},P^0) \rho({\rm d}\theta) + \mathcal{K}(\rho,\pi) \bigg\}
 \\ \label{defELBO} 
 & = \argmax_{\rho \in \mathcal{F}} \bigg\{ \alpha \int \ell_n(\theta) \rho(d\theta) - \mathcal{K}\big(\rho,\pi\big) \bigg\}
 \\ \label{defpasELBO} 
 & = \argmin_{\rho \in \mathcal{F}} \bigg\{ - \alpha \sum_{i=1}^n \int \log\bigg({\sum_{j=1}^K p_j q_{\theta_j}(X_i)}\bigg) \rho(d\theta) + \mathcal{K}\big(\rho_p,\pi_p\big) + \sum_{j=1}^K \mathcal{K}\big(\rho_j,\pi_j\big) \bigg\}.
\end{align}
The quantity maximized in~\eqref{defELBO} is called the ELBO in the litterature (ELBO stands for Evidence Lower Bound), and many authors actually take this as the definition of VB \cite{Review}.

\begin{rmk}
In practice, the choice of $\alpha$ is not staightforward. Depending on the objective, some heuristic might be available: for example,~\cite{syring2015scaling} proposed a nice method to calibrate $\alpha$ in order to get confidence intervals on a parameter of interest. More generally, cross-validation can give good results. We have to acknowledge that there is no universal method to calibrate $\alpha$. This could lead the reader to the idea that the proper Bayesian approach ($\alpha=1$) is simpler to use. We insist on the fact $\alpha=1$ can produce catastrophic results in case of misspecification~\cite{grunwaldmisspecifiation}, while we present below some results in the misspecified case with $\alpha<1$.  We believe that the calibration of $\alpha$ is a very important research direction.
\end{rmk}

\vspace{0.2cm}

\section{Variational Bayes estimation of a mixture}
\label{sec:main-res}

\vspace{0.2cm}
\subsection{A PAC-Bayesian inequality}
\label{subsec:pacbayes}
\vspace{0.2cm}

We start with a result for general mixtures. Later in this section we provide corollaries obtained by applying this theorem to special cases: mixture of multinomials and Gaussian mixtures.

\begin{thm}
 \label{thm-general}
 For any $\alpha\in(0,1)$,
 
 \begin{align*}
 & \mathbb{E} \bigg[ \int D_{\alpha}\bigg( \sum_{j=1}^K p_j q_{\theta_j} ,\sum_{j=1}^K p_j^0 q_{\theta_j^0} \bigg) \tilde{\pi}_{n,\alpha}(d\theta|X_1^n) \bigg]
 \\
 & \quad \leq \inf_{\rho_p \bigotimes_{j=1}^{K} \rho_j \in \mathcal{F}} \bigg\{ \frac{\alpha}{1-\alpha} \bigg[ \int \mathcal{K}(p^0,p) \rho_p(dp) + \sum_{j=1}^K \int \mathcal{K}(q_{\theta_j^0},q_{\theta_j}) \rho_j(d\theta_j) \bigg] + \frac{\mathcal{K}(\rho_p,\pi_p) + \sum_{j=1}^K \mathcal{K}(\rho_j,\pi_j)}{n(1-\alpha)}  \bigg\}.
 \end{align*}
 
 \noindent
As a special case, when there exists $r_{n,K}$ such that there is are distributions $\rho_{p,n} \in \mathcal{M}_1^+(\mathcal{S}_K)$ and $\rho_{j,n} \in \mathcal{M}_1^+(\Theta)$ ($j=1,...,K$) such that for $j=1,...,K$
\begin{equation}
  \int \mathcal{K}(p^0,p) \rho_{p,n}(dp) \leq K r_{n,K} , \hspace{0.2cm}
  \int \mathcal{K}(q_{\theta_j^0},q_{\theta_j}) \rho_{j,n}(d\theta_j) \leq r_{n,K}
  \label{cond:1}
\end{equation}
and 
\begin{equation}
  \mathcal{K}(\rho_{p,n},\pi_p) \leq K n r_{n,K} , \hspace{0.2cm}
  \mathcal{K}(\rho_{j,n},\pi_j) \leq n r_{n,K},
  \label{cond:2}
\end{equation}
then for any $\alpha \in (0,1)$

$$
\mathbb{E} \bigg[ \int D_{\alpha}\bigg( \sum_{j=1}^K p_j q_{\theta_j} ,\sum_{j=1}^K p_j^0 q_{\theta_j^0} \bigg) \tilde{\pi}_{n,\alpha}(d\theta|X_1^n) \bigg] \leq \frac{1+\alpha}{1-\alpha} 2K r_{n,K}.
$$

\end{thm}

\vspace{0.2cm}
The proof is given in Section \ref{sec:proofs}. This theorem provides the consistency of the Variational Bayes for mixture models as soon as \eqref{cond:1} and \eqref{cond:2} are satisfied. In \cite{Tempered}, the authors use similar conditions ((3) and (4) in their Theorem 2.6), and show that they are strongly linked to the assumptions on the prior used by \cite{ghosal2000convergence,bhattacharya2016bayesian} to derive concentration of the posterior. Thus they cannot be removed in general. Theorem \ref{thm-general} states that finding $r_{n,K}$ fulfilling \eqref{cond:1} and \eqref{cond:2} independently for the weights and for each component is sufficient to obtain the rate of convergence $Kr_{n,K}$ of the VB estimator towards the true distribution.

\vspace{0.2cm}
Clearly, the theorem cannot be directly extended to the case $\alpha=1$. As discussed above, the case $\alpha=1$ is studied in~\cite{Chicago}: it requires a testing condition in addition to the prior mass condition given by~\eqref{cond:1} and~\eqref{cond:2}. In the case of mixtures, such a testing condition was studied in~\cite{kruijer2010adaptive}.

\vspace{0.2cm}
Note that there always exists a distribution $\rho_{p,n} \in \mathcal{M}_1^+(\mathcal{S}_K)$ such that the two quantities corresponding to the weights $\int \mathcal{K}(p^0,p) \rho_{p,n}(dp)$ and $\mathcal{K}(\rho_{p,n},\pi_p)$ are bounded as required in Theorem \ref{thm-general} for $r_{n,K}=\frac{4\log(nK)}{n}$ when the chosen prior is a Dirichlet distribution $\pi_p = \mathcal{D}_K(\alpha_1,...,\alpha_K)$ under some minor restriction on $\alpha_1,...,\alpha_K$. This result summarized below for any $K\geq2$ helps find explicit rates of convergence for the VB approximation.

\begin{lemma}
\label{lemma-Dirichlet}

For $r_{n,K}= \frac{4\log(nK)}{n}$ and a prior $\pi_p = \mathcal{D}_K(\alpha_1,...,\alpha_K) \in \mathcal{M}_1^+(\mathcal{S}_K)$ with $\frac{2}{K} \leq \alpha_j \leq 1$ for $j=1,...,K$, we can find a distribution $\rho_{p,n} \in \mathcal{M}_1^+(\mathcal{S}_K)$ such that

\begin{equation*}
  \int \mathcal{K}(p^0,p) \rho_{p,n}(dp) \leq K r_{n,K}
\end{equation*}
and 
\begin{equation*}
  \mathcal{K}(\rho_{p,n},\pi_p) \leq K n r_{n,K}.
\end{equation*}

\end{lemma}

\vspace{0.2cm}
Thus, conditions \eqref{cond:1} and \eqref{cond:2} concerning the mixture components are always satisfied for guaranteeing consistency and obtaining convergence rates of the Variational Bayes procedure.

\vspace{0.2cm}
\begin{rmk}
When $K=1$, Lemma \ref{lemma-Dirichlet} does not apply as $\frac{2}{K}>1$. Nevertheless, as there is only one component, then any $p \in \mathcal{S}_K$ is equal to $1$ and the two conditions are immediately satisfied for any prior $\pi_p$ and any rate $r_{n,K}$ with $\rho_{p,n}=\pi_p$.
\end{rmk}

\vspace{0.2cm}
The central idea of the proof of Lemma \ref{lemma-Dirichlet} (given in details in Section \ref{sec:proofs}) is to consider the ball $\mathcal{B}$ centered at $p^0$ of radius $Kr_{n,K}$ defined as:
$$
\mathcal{B} = \bigg\{p \in \mathcal{S}_K / \hspace{0.2cm}  \mathcal{K}(p^0,p) \leq K r_{n,K}\bigg\}.
$$
Hence, when considering the restriction $\rho_{p,n} \in \mathcal{M}_1^+(\mathcal{S}_K)$ of $\pi_p$ to $\mathcal{B}$, condition \eqref{cond:1} is trivially satisfied and condition \eqref{cond:2} is restricted to 
$$
\pi_p(\mathcal{B}) \geq e^{-nKr_{n,K}}.
$$
This is a very classical assumption stated in many papers to study the concentration of the posterior \cite{ghosal2000convergence,Tempered,Chicago}. However, the computation of such a prior mass $\pi_p(\mathcal{B})$ is a major difficulty. Lemma 6.1 in \cite{ghosal2000convergence} treated the case of $L_1$-balls for Dirichlet priors. Since then, only a few papers in the literature addressed this issue. Our result extends the work in \cite{ghosal2000convergence} to KL-balls, which is of great interest in our study. Moreover, the range of Dirichlet priors for which Lemma \ref{lemma-Dirichlet} is applicable is the same as the one in \cite{ghosal2000convergence}.

\vspace{0.2cm}
We conclude Subsection~\ref{subsec:pacbayes} by a short discussion on the implementation of the VB approximation. Indeed, VB methods are meant to be practical objects, so there would be no point in proving the consistency of a VB approximation that would not be computable in practice. Many algorithms have been studied in the literature, with good performances -- see \cite{Review} and the references therein. In the case of mean-field approximation, the most popular method is to optimize iteratively with respect to all the independent components. Here this might seem difficult: it is indeed as difficult as maximizing the likelihood of a mixture. But a trick widely used in practice (see for example Section 7 in \cite{DoLemmaSignalProcessing}) is to use the equality
$$ \textnormal{for any } i=1,...,K \hspace{0.7cm}
- \log\bigg(\sum_{j=1}^K p_j q_{\theta_j}(X_i)\bigg) = \min_{\omega^i \in \mathcal{S}_K} \bigg\{ - \sum_{j=1}^K \omega^i_j \log(p_j q_{\theta_j}(X_i)) + \sum_{j=1}^K \omega^i_j \log(\omega^i_j) \bigg\}.
$$
This equality is once again a consequence of Lemma~\ref{thm-dv} (take $\textbf{E}=\{1,...,K\}$, $\lambda=(1/K,...,1/K)$ and  $h(j)=\log(p_jq_{\theta_j}(X_i))$).
This leads to the program:
\begin{equation*}
\begin{split}
\min_{\rho \in \mathcal{F}, \hspace{0.1cm} w \in \mathcal{S}_K^n} \bigg\{ - \alpha \sum\limits_{i=1}^n \sum\limits_{j=1}^K \omega_j^i & \bigg( \int \log(p_j) \rho_p(dp) + \int \log(q_{\theta_j}(X_i)) \rho_j(d\theta_j) \bigg) \\
& + \alpha \sum\limits_{i=1}^n \sum\limits_{j=1}^K \omega_j^i \log(\omega_j^i)
+ \mathcal{K}(\rho_p,\pi_p) + \sum\limits_{j=1}^K \mathcal{K}(\rho_j,\pi_j) \bigg\}.
\end{split} 
\end{equation*}
This version can be solved by coordinate descent, see Algorithm 1. Update formulas once again follow from Lemma \ref{thm-dv} (for instance, line 7 can be obtained with $\textbf{E}=\{1,...,K\}$, $\lambda=(1/K,...,1/K)$ and  $h(j)=\int \log(p_j) \rho_p(dp) + \int \log(q_{\theta_j}(X_i)) \rho_j(d\theta_j)$, more details are provided in Section~\ref{sec:proofs}). This algorithm is, in the case $\alpha=1$, exactly equivalent to the popular CAVI algorithm \cite{Plage,Review,hoffman2013stochastic}, where the $\omega^i_j$'s are interpreted as the posterior means of the latent variables $Z^i_j$'s. A very short numerical study is provided in the Supplementary Material but note that CAVI has already been extensively tested in practice \cite{Review}.

\begin{algorithm}
\caption{Coordinate Descent Variational Bayes for mixtures}
\begin{algorithmic}[1]
\State \textbf{Input}: a dataset $(X_1,...,X_n)$, priors $\pi_p$,$\{\pi_{j}\}_{j=1}^K$ and a family  $\{q_\theta/\theta \in \Theta\}$
\State \textbf{Output}: a variational approximation $\rho_p(p)\prod_{j=1}^{K}\rho_{j}(\theta_j)$
\State \textbf{Initialize} variational factors $\rho_p$, $\{\rho_{j}\}_{j=1}^K$
\State \textbf{until} convergence of the objective function \textbf{do}
   \For {$i=1,...,n$}
       \For {$j=1,...,K$}
          \State set $w^i_j=\exp\bigg( \int \log(p_j) \rho_p(dp) + \int \log(q_{\theta_j}(X_i)) \rho_j(d\theta_j) \bigg)$
      \EndFor
   \State normalize $(w^i_j)_{1\leq j\leq K}$
   \EndFor
   \State
   set $\rho_p(dp) \propto \exp\bigg( \alpha \sum\limits_{i=1}^n \sum\limits_{j=1}^K \omega^i_j \log(p_j) \bigg) \pi_p(dp)$
   \For {$j=1,...,K$}
      \State set $\rho_j(d\theta_j) \propto \exp\bigg( \alpha \sum\limits_{i=1}^n \omega^i_j \log(q_{\theta_j}(X_i)) \bigg) \pi_j(d\theta_j)$
   \EndFor
\end{algorithmic}
\end{algorithm} 

\vspace{0.2cm}
\subsection{Application to multinomial mixture models}
\label{sec:main-res:2}
\vspace{0.2cm}

We present in this section an application to the multinomial mixture model frequently used for text clustering \cite{multinomialmixturemodel}, transport schedule analysis \cite{carel2017simultaneous} and others. The parameter space is the $V-1$ dimensional simplex $\Theta=\mathcal{S}_V$ with $V$ a positive integer, and $q_{\theta}(X)=\prod_{v=1}^V \theta_{vj}^{\mathds{1}(X=v)}$ for any $\theta \in \Theta$. We choose conjugate Dirichlet priors as in \cite{multinomialmixturemodel} $\pi_p = \mathcal{D}_K(\alpha_1,...,\alpha_K)$ and $\pi_j = \mathcal{D}_V(\beta_{1},...,\beta_{V})$ with $\frac{2}{K} \leq \alpha_j \leq 1$ for $j=1,...,K$ and $\frac{2}{V} \leq \beta_{\ell} \leq 1$ for $\ell=1,...,V$.

\vspace{0.2cm}
The following corollary of Theorem \ref{thm-general} states that convergence of the VB approximation for the multinomial mixture model is achieved at rate $\frac{KV\log(nV)}{n}$ as soon as $V^{V} \geq K$, which is the case in many text mining models such as Latent Dirichlet Allocation \cite{LDA} for which the size of the vocabulary is very large:

\begin{cor} 
\label{thm-multinomial}
For any $\alpha \in (0,1)$,
\begin{equation*}
   \begin{split}
       \mathbb{E} \bigg[ \int D_{\alpha}\bigg( \sum_{j=1}^K p_j q_{\theta_j} ,\sum_{j=1}^K & p_j^0 q_{\theta_j^0} \bigg) \tilde{\pi}_{n,\alpha}(d\theta|X_1^n) \bigg] \leq \frac{1+\alpha}{1-\alpha} \bigg[ \frac{8KV\log(nV)}{n} \bigvee \frac{8K\log(nK)}{n} \bigg]
   \end{split}. 
\end{equation*}
\end{cor}

\vspace{0.2cm}
The proof is in Section \ref{sec:proofs}. We also specialize Algorithm 1 to the present setting (see Algorithm 2). Here $\psi$ denotes the Digamma function, $\psi(x)=\frac{{\rm d}}{{\rm d}x} \log[\Gamma(x)] $ where $\Gamma$ stands for the Gamma function $\Gamma(x)=\int_{0}^{\infty} \exp(-t) t^{x-1} {\rm d}t$.

\begin{algorithm}
\caption{Coordinate Descent Variational Bayes for multinomial mixtures}
\begin{algorithmic}[1]
\State \textbf{Initialize} variational parameters $(\phi_1,...,\phi_K) \in \mathbb{R}_+^K$, $(\gamma_{1j},...,\gamma_{Vj}) \in \mathbb{R}_+^V$ and corresponding variational distributions $\rho_p=\mathcal{D}_K(\phi_1,...,\phi_K)$, $\rho_{j}=\mathcal{D}_V(\gamma_{1j},...,\gamma_{Vj})$ for $j=1,...,K$
\State \textbf{until} convergence of the objective function \textbf{do}
   \For {$i=1,...,n$}
       \For {$j=1,...,K$}
          \State set $w^i_j=\exp\bigg( \psi(\phi_j)-\psi(\sum\limits_{\ell=1}^K \phi_{\ell}) +  \psi(\gamma_{X_i,j}) - \psi\big(\sum\limits_{v=1}^V
          \gamma_{vj}\big) \bigg)$
      \EndFor
   \State normalize $(w^i_j)_{1\leq j\leq K}$
   \EndFor
   \State set $\phi_j = \alpha_j+\alpha \sum\limits_{i=1}^n\omega_{j}^i$ \hspace{0.1cm} \textnormal{for} \hspace{0.1cm} $j=1,...,K$
   \State
   set $\rho_{p}=\mathcal{D}_K(\phi_1,...,\phi_K)$   
   \For {$j=1,...,K$}
      \State set $\gamma_{vj}=\beta_v+\alpha \sum\limits_{i=1}^n \omega_{j}^i \mathds{1}(X_i=v)$ \hspace{0.1cm} \textnormal{for} \hspace{0.1cm} $v=1,...,V$
      \State set $\rho_{j}=\mathcal{D}_V(\gamma_{1j},...,\gamma_{Vj})$
   \EndFor
\end{algorithmic}
\end{algorithm}

\vspace{0.2cm}
\subsection{Application to Gaussian mixture models}
\label{sec:main-res:3}
\vspace{0.2cm}

Let us now address the case of the Gaussian mixture model. This is one of the most popular mixture models for many applications including model based clustering \cite{Bouveyron2014,mcnicholas2016} and VB approximations have been studied in depth for this model \cite{VariationallearningforGaussianmixturemodels}. First, we will give rates of convergence of the VB approximation of the tempered posterior when the variance is known, and then when the variance is unknown. 

\vspace{0.2cm}

First, we consider mixtures of $V^2$-variance Gaussians. The parameter space is $\Theta=\mathbb{R}$, and each component $j=1,...,K$ is parameterized by its mean $\theta_j=\mu_j$. We  select priors $\pi_p = \mathcal{D}_K(\alpha_1,...,\alpha_K)$ and $\pi_j = \mathcal{N}(0,\mathcal{V}^2)$ with $\frac{2}{K} \leq \alpha_j \leq 1$ for $j=1,...,K$ and $\mathcal{V}^2 > 0$. The following result gives a rate of convergence $Kr_{n,K}$ of the VB approximation:

\begin{cor}
\label{thm-Gaussian}
Let us define 
$r_{n,K} = \frac{4\log(nK)}{n} \bigvee_{j=1}^{K} \frac{1}{n} \bigg[ \frac{1}{2} \log\bigg(\frac{n}{2}\bigg) + \frac{V^2}{n\mathcal{V}^2} + \log\bigg(\frac{\mathcal{V}}{V}\bigg) + \frac{(\mu_j^0)^2}{2\mathcal{V}^2} - \frac{1}{2} \bigg] $. Then, for any $\alpha \in (0,1)$,
$$
\mathbb{E} \bigg[ \int D_{\alpha}\bigg( \sum_{j=1}^K p_j q_{\theta_j} ,\sum_{j=1}^K p_j^0 q_{\theta_j^0} \bigg) \tilde{\pi}_{n,\alpha}(d\theta|X_1^n) \bigg] \leq \frac{1+\alpha}{1-\alpha} 2Kr_{n,K}.
$$

\end{cor}

One can see that for $n$ large enough, the convergence rate is $\frac{K\log(nK)}{n}$, which comes from the estimation of the weights of the mixture.

\vspace{0.2cm}

We can also provide a similar result when the variance of each component is unknown. The convergence rate remains the same, and is entirely characterized by the weights consistency rate. The parameter space is now $\Theta=\mathbb{R} \times (0,+\infty)$, and each component $j=1,...,K$ is parameterized by its pair mean/variance $\theta_j=(\mu_j,\sigma_j^2)$. We consider again a Dirichlet prior $\pi_p = \mathcal{D}_K(\alpha_1,...,\alpha_K)$ with $\frac{2}{K} \leq \alpha_j \leq 1$ for $j=1,...,K$ on $p\in \mathcal{S}_K$, and we will provide our results for two different priors $\pi_j$ frequently used in the literature: a Normal-Inverse-Gamma prior \cite{NormalInverseGammaPriorForMixtureModels} and a factorized prior \cite{Watier99usinggaussian}. We define the Normal-Inverse-Gamma distribution as follows:
\begin{dfn}
 The Normal-Inverse-Gamma $\mathcal{NIG}(\mu,\theta^2,a,b)$ is the distribution which density $w$ with respect to Lebesgue measure is defined by $w(x,y)=g(x|\mu,\frac{y}{\theta^2}) h(y|a,b)$, where $g(.|\mu,\sigma^2)$ is the density function of a Gaussian distribution of mean $\mu$ and variance $\sigma^2$, and $h(.|a,b)$ is the density distribution of an Inverse-Gamma of parameters $a$ and $b$.
\end{dfn}

\begin{cor}
\label{thm-Gaussian-variance}
Let us fix $\alpha \in (0,1)$. 
\begin{itemize}
    \item For a Normal-Inverse-Gamma prior $\pi_j = \mathcal{NIG}(0,\mathcal{V}^{-2},1,\gamma^2)$ for each $j=1,...,K$. With \newline $r_{n,K} =  \frac{4\log(nK)}{n} \bigvee_{j=1}^{K} \frac{1}{n} \bigg[  2 \log(n\sqrt{\mathcal{V}}) + \frac{1}{2n\mathcal{V}^2} + \frac{(\mu_j^0)^2}{2(\sigma_j^0)^2\mathcal{V}^2} + \log\bigg(\frac{(\sigma_j^0)^2}{\gamma^2}\bigg) + \frac{\gamma^2}{(\sigma_j^0)^2}  - \frac{1}{2}\log(2\pi)  \bigg]$,
    $$
\mathbb{E} \bigg[ \int D_{\alpha}\bigg( \sum_{j=1}^K p_j q_{\theta_j} ,\sum_{j=1}^K p_j^0 q_{\theta_j^0} \bigg) \tilde{\pi}_{n,\alpha}(d\theta|X_1^n) \bigg] \leq \frac{1+\alpha}{1-\alpha} 2Kr_{n,K}.
$$
    \item For the factorized prior $\pi_j = \mathcal{N}(0,\mathcal{V}^2) \bigotimes \mathcal{IG}(1,\gamma^2)$ for each $j=1,...,K$. With \\ $r_{n,K} = \frac{4\log(nK)}{n} \bigvee_{j=1}^{K} \frac{1}{n} \bigg[ 2 \log(n\sqrt{\mathcal{V}}) + \frac{(\sigma_j^0)^2}{2n\mathcal{V}^2} + \frac{(\mu_j^0)^2}{2\mathcal{V}^2} + \frac{1}{2} \log\bigg(\frac{(\sigma_j^0)^2}{\gamma^4}\bigg) + \frac{\gamma^2}{(\sigma_j^0)^2}  - \frac{1}{2}\log(2\pi) \bigg] $,
    $$
\mathbb{E} \bigg[ \int D_{\alpha}\bigg( \sum_{j=1}^K p_j q_{\theta_j} ,\sum_{j=1}^K p_j^0 q_{\theta_j^0} \bigg) \tilde{\pi}_{n,\alpha}(d\theta|X_1^n) \bigg] \leq \frac{1+\alpha}{1-\alpha} 2Kr_{n,K}.
$$
\end{itemize}

\end{cor}

One can see that even when the variance has to be estimated, the convergence rate still achieves $\frac{K\log(nK)}{n}$ for $n$ large enough, whatever the form of the prior - factorized or not.

\vspace{0.2cm}
We give in Algorithm 3 a version of Algorithm 1 for unit-variance Gaussian mixtures with priors $\pi_p = \mathcal{D}_K(\alpha_1,...,\alpha_K)$ and $\pi_j = \mathcal{N}(0,\mathcal{V}^2)$ where $\frac{2}{K} \leq \alpha_j \leq 1$ for $j=1,...,K$ and $\mathcal{V}^2 > 0$.

\begin{algorithm}
\caption{Coordinate Descent Variational Bayes for unit-variance Gaussian mixtures}
\begin{algorithmic}[1]
\State \textbf{Initialize} variational parameters $(\phi_1,...,\phi_K) \in (\mathbb{R}_+^*)^K$, $(n_j,s_j^2) \in \mathbb{R} \times \mathbb{R}_+^*$ and corresponding variational distributions $\rho_p=\mathcal{D}_K(\phi_1,...,\phi_K)$, $\rho_{j}=\mathcal{N}(n_j,s_j^2)$ for $j=1,...,K$
\State \textbf{until} convergence of the objective function \textbf{do}
   \For {$i=1,...,n$}
       \For {$j=1,...,K$}
          \State set $w^i_j=\exp\bigg( \psi(\phi_j)-\psi(\sum\limits_{\ell=1}^K \phi_{\ell}) - \frac{1}{2}\big\{ s_j^2 + (n_j-X_i)^2 \big\} \bigg)$
      \EndFor
   \State normalize $(w^i_j)_{1\leq j\leq K}$
   \EndFor
   \State
   set $\phi_j = \alpha_j+\alpha \sum\limits_{i=1}^n\omega_{j}^i$ \hspace{0.1cm} \textnormal{for} \hspace{0.1cm} $j=1,...,K$
   \State
   set $\rho_{p}=\mathcal{D}_K(\phi_1,...,\phi_K)$
   \For {$j=1,...,K$}
      \State set $n_j=\frac{\alpha \sum_{i=1}^n \omega_{j}^iX_i}{1/\mathcal{V}^2+\alpha \sum_{i=1}^n \omega_{j}^i}$ and $s_j^2=\frac{1}{1/\mathcal{V}^2+\alpha \sum_{i=1}^n \omega_{j}^i}$
      \State set $\rho_{\mu,j}=\mathcal{N}(n_j,s_j^2)$
   \EndFor
\end{algorithmic}
\end{algorithm}

\vspace{0.2cm}
\subsection{Extension to the misspecified case}
\label{sec:main-res:4}
\vspace{0.2cm}

From now we do not assume any longer that the true distribution $P^0$ belongs to the $K$-mixtures model. We still consider a prior $\pi=\pi_p \bigotimes_{j=1}^{K}\pi_j$ on $\theta \in \Theta_K$ for which $\pi_p \in \mathcal{M}_1^+(\mathcal{S}_K)$ and $\pi_j \in \mathcal{M}_1^+(\Theta)$ for $j=1,...,K$.

\vspace{0.2cm}
For some value $r_{n,K}$, we introduce the set $\Theta_K(r_{n,K})$ of parameters $\theta^* \in \Theta_K$ such that:
\begin{itemize}
    \item there exists a set $\mathcal{A}_{n,K} \subset \mathcal{S}_K$ satisfying:
        \begin{itemize}
            \item for each $p\in \mathcal{A}_{n,K}$, for each $j=1,...,K$, $\hspace{0.2cm}$ $\log\big(\frac{p_j^*}{p_j}\big) \leq Kr_{n,K}$,
            \item $\pi_p(\mathcal{A}_{n,K}) \geq e^{-nKr_{n,K}}$.
        \end{itemize}
    \item there are distributions $\rho_{j,n} \in \mathcal{M}_1^+(\Theta)$ ($j=1,...,K$) such that for $j=1,...,K$:
\begin{equation}
  \int \mathbb{E}\bigg[\log\bigg(\frac{q_{\theta_{j}^*}(X)}{q_{\theta_{j}}(X)} \bigg)\bigg]  \rho_{j,n}(d\theta_{j}) \leq r_{n,K} \hspace{0.2cm},  \hspace{0.5cm}
  \mathcal{K}(\rho_{j,n},\pi_{j}) \leq n r_{n,K}.
  \label{cond:3}
\end{equation}
\end{itemize}

\vspace{0.2cm}
Let us discuss this definition. To begin with, the first item of the definition of $\Theta_K(r_{n,K})$ can seem quite restrictive. It is even a much more stronger assumption than \eqref{cond:1} and \eqref{cond:2}. Nevertheless, the way to find the required measures $\rho_{p,n}$ in Lemma \ref{lemma-Dirichlet} in the well-specified case implies constructing in the proof such sets $\mathcal{A}_{n,K}$ for the true weight parameter $p^0$. As a consequence, it might seem reasonable to replace conditions \eqref{cond:1} and \eqref{cond:2} by the first part of the definition of $\Theta_K(r_{n,K})$. On the other hand, the condition given by \eqref{cond:3} looks like those of Theorem 2.7 in \cite{Tempered}. Once again, the difference is that inequalities must be satisfied here for each component. A condition on both the true distribution $P^0$ and the parameter $\theta^*$ considered is required through the expectation term. Besides, condition \eqref{cond:3} is equivalent to \eqref{cond:1} and \eqref{cond:2} when the model is well-specified.

\vspace{0.2cm}

\begin{thm}
\label{thm-misspecified}
For any $\alpha \in (0,1)$,
$$
\mathbb{E} \bigg[ \int D_{\alpha}\bigg( \sum_{j=1}^K p_j Q_{\theta_j} ,P^0 \bigg) \tilde{\pi}_{n,\alpha}(d\theta|X_1^n) \bigg] \leq \frac{\alpha}{1-\alpha}  \inf_{ \theta^* \in \Theta_K(r_{n,K}) } \mathcal{K}(P^0,P_{\theta^*}) + \frac{1+\alpha}{1-\alpha} 2K r_{n,K}.
$$

\end{thm}

\begin{rmk}
If there is no $r_{n,K}$ such that $\Theta_K(r_{n,K})$ is not empty, then the right-hand side is equal to infinity (by convention) for any value of $r_{n,K}$ and the inequality is useless. Nevertheless, this is not the case in models used in practice. We show an example below.
\end{rmk}

It is worth mentioning that even if this is not exactly an oracle inequality as the risk function in the left-hand side ($\alpha$-Renyi divergence) is lower than the right-hand side one (Kullback-Leibler divergence), but the theorem still remains of great interest. Indeed, when the minimizer of $\mathcal{K}(P^0,P_{\theta})$ with respect to $\theta \in \Theta_K(r_{n,K})$ exists and is such that the corresponding Kullback-Leibler divergence is small, then our oracle inequality is informative as it gives a small bound on the expected risk of the Variational Bayes.

\vspace{0.2cm}
To illustrate the relevance of Theorem \ref{thm-misspecified}, we provide the following result that is applicable for a wide range of generating distributions when considering the family of unit-variance Gaussian mixtures with priors $\pi_p = \mathcal{D}_K(\alpha_1,...,\alpha_K) \in \mathcal{M}_1^+(\mathcal{S}_K)$ with $\frac{2}{K} \leq \alpha_j \leq 1$ for $j=1,...,K$ ($K\geq2)$ and $\pi_j = \mathcal{N}(0,\mathcal{V}^2) \in \mathcal{M}_1^+(\mathbb{R})$ for $j=1,...,K$ with $\mathcal{V}^2 > 0$:

\begin{cor}
\label{thm-misspecified-Gaussian}
Assume that the true distribution $P^0$ is such that $\mathbb{E}|X| < +\infty$. Let $L>0$. \\
For $r_{n,K}=\frac{4\log(nK)}{n}  \bigvee_{j=1}^{K} \frac{1}{n} \bigg[ \frac{1}{2} \log\bigg(\frac{n}{2}\bigg) + \frac{1}{n\mathcal{V}^2} + \log\big({\mathcal{V}}\big) + \frac{L^2}{2\mathcal{V}^2} - \frac{1}{2} \bigg]$, we get $\mathcal{S}_K \times [-L,L]^K \subset \Theta_K(r_{n,K})$  and for any $\alpha \in (0,1)$,
$$
\mathbb{E} \bigg[ \int D_{\alpha}\bigg( \sum_{j=1}^K p_j Q_{\theta_j} ,P^0 \bigg) \tilde{\pi}_{n,\alpha}(d\theta|X_1^n) \bigg] \leq \frac{\alpha}{1-\alpha} \inf_{ \theta^* \in \mathcal{S}_K \times [-L,L]^K } \mathcal{K}(P^0,P_{\theta^*}) + \frac{1+\alpha}{1-\alpha} 2K r_{n,K}.
$$
\end{cor}

\begin{rmk}
If the true distribution is a mixture of unit-variance Gaussians with components means between $-L$ and $L$, then $\mathbb{E}|X| < +\infty$ and the first term of the right-hand side of the inequality is equal to zero, which gives directly for any $\alpha \in (0,1)$,
$$
\mathbb{E} \bigg[ \int D_{\alpha}\bigg( \sum_{j=1}^K p_j Q_{\theta_j} ,P^0 \bigg) \tilde{\pi}_{n,\alpha}(d\theta|X_1^n) \bigg] \leq \frac{1+\alpha}{1-\alpha} 2K r_{n,K}.
$$
\end{rmk}

\vspace{0.2cm}
\section{Variational Bayes model selection}
\label{sec:selection}
\vspace{0.2cm}

In this section, we extend the problem to a larger family of distributions. We want to model the generating distribution $P^0$ using mixtures with an unknown number of components in a possibly misspecified setting. Thus, we consider a countable collection $ \left\lbrace \mathcal{M}_K /  K \in \mathbb{N}^* \right\rbrace$  of statistical mixture models 
$$ 
\mathcal{M}_K = \left\lbrace P_{\theta_K} = \sum_{j=1}^K p_{j,K} Q_{\theta_{j,K}} \hspace{0.2cm} / \hspace{0.2cm} \theta_{K} \in \Theta_K \right\rbrace
$$ 
with $\Theta_K=\mathcal{S}_K \times \Theta^K$, $\mathcal{S}_K=\{ p_K=(p_{1,K},...,p_{K,K}) \in [0,1]^K / \sum_{j=1}^K p_{j,K} = 1 \}$ and the general notation $\theta_K=(p_{K},\theta_{1,K},...,\theta_{K,K})$. We would like to emphasize that the notations are slightly different
 as the size of each component parameter depends on the model complexity $K$. The entire parameter space $\Omega$ is the union of all parameter spaces $\Theta_K$ associated with each model index $K$: $ \Omega = \cup_{K=1}^{\infty} \Theta_K$, and we can think of a whole statistical model $\mathcal{M}=\cup_{K=1}^\infty \mathcal{M}_K$ as the union of all collections $\mathcal{M}_K$. First, we can notice that different models $\mathcal{M}_K$ never overlap as parameters in each one do not have the same length. Nonetheless, parameters in complex models (models $\mathcal{M}_K$ with large $K$) can be sparse and therefore contain the "same information" as parameters in less complex ones, i.e. can lead to the same distribution $P_{\theta}$.

\vspace{0.2cm}
The prior specification is a crucial point. As mentioned above, each parameter depends on the number of components. Then, we specify a prior weight $\pi_K$ assigned to the model $\mathcal{M}_K$ and a conditional prior $\Pi_K(.)$ on $\theta_K \in \Theta_K$ given model $\mathcal{M}_K$. More precisely, we define our conditional prior on $\theta_K=(p_{K},\theta_{1,K},...,\theta_{K,K})$ as follows: given $K$, the weight parameter $p_K=(p_{1,K},...,p_{K,K})$ is supposed to follow a distribution $\pi_{p,K}$ on $\mathcal{M}_1^+(\mathcal{S}_K)$; finally, given $K$, we set independent priors $\pi_{j,K}$ for the component parameters $\theta_{j,K}$ where each $\pi_{j,K}$ is a probability distribution on $\mathcal{M}_1^+(\Theta)$. In a nutshell:

$$ \pi = \sum_{K=1}^{+\infty} \pi_K \Pi_K $$
with
$$\Pi_K(\theta_K)= \pi_{p,K}(p_K) \prod_{j=1}^{K}\pi_{j,K}(\theta_{j,K}).$$

\vspace{0.2cm}
We have to adapt the notations for the VB approximations. The tempered posteriors $\pi_{n,\alpha}^K(.|X_1^n)$ on parameter $\theta_K \in \Theta_K$ given model $\mathcal{M}_K$, is defined again as
\[
\pi_{n,\alpha}^K(d\theta_K|X_1^n) \propto L_{n}(\theta_K)^\alpha \Pi_K(d\theta_K).
\]
The Variational Bayes $\tilde{\pi}_{n,\alpha}^K(.|X_1^n)$ is the projection of the tempered posterior onto some set $\mathcal{F}_K$ following the mean-field assumption: the variational factor corresponding to the weight parameter $p_K=(p_{1,K},...,p_{K,K})$ is any distribution $\rho_{p}$ on $\mathcal{M}_1^+(\mathcal{S}_K)$; besides, we consider independent variational distributions $\rho_{j}(\theta_{j,K})$ for the component parameters $\theta_{j,K}$ where each $\rho_{j}$ is a probability distribution on $\mathcal{M}_1^+(\Theta)$. Then, $\mathcal{F}_K=\{\rho_p \bigotimes_{j=1}^{K} \rho_j / \rho_p \in \mathcal{M}_1^+(\mathcal{S}_K), \hspace{0.1cm} \rho_j \in \mathcal{M}_1^+(\Theta) \hspace{0.1cm} \forall j=1,...,K \}$, and
\[
\tilde{\pi}_{n,\alpha}^K(.|X_1^n) = \argmin_{\rho_K \in \mathcal{F}_K} \mathcal{K}\bigg(\rho_K,\pi_{n,\alpha}^K(.|X_1^n)\bigg).
\]
We recall that an alternative way to define the variational estimate is to use the Evidence Lower Bound via the optimization program (\ref{defELBO}):
\[
\tilde{\pi}_{n,\alpha}^K(.|X_1^n) = \argmax_{\rho_K \in \mathcal{F}_K} \bigg\{ \alpha \int \ell_n(\theta_K) \rho_K(d\theta_K) - \mathcal{K}\big(\rho_K,\Pi_K\big) \bigg\}
\]
where the function inside the argmax operator is the ELBO $\mathcal{L}(\rho_K)$. For simplicity, we will just call ELBO $\mathcal{L}(K)$ the closest approximation to the log-evidence, i.e. the value of the lower bound evaluated in its maximum:
$$
\mathcal{L}(K)=\alpha \int \ell_n({\theta}_K) \tilde{\pi}_{n,\alpha}^K(d\theta_K|X_1^n)-\mathcal{K}(\tilde{\pi}_{n,\alpha}^K(.|X_1^n),\Pi_K).
$$

\vspace{0.2cm}
The objective is to propose a data-driven estimate $\hat{K}$ of the number of components from which we will pick up our final VB estimate $\tilde{\pi}_{n,\alpha}^{\hat{K}}(.|X_1^n)$ and derive an oracle inequality in the spirit of \cite{Massart}. It is stated in \cite{Review} that $\argmax_{K\geq 1}  \mathcal{L}(K)$ is widely used in practice, without any theoretical justification. We propose
$$
\hat{K} = \argmax_{K\geq 1} \bigg\{ \mathcal{L}(K) - {\log\bigg(\frac{1}{\pi_K}\bigg)} \bigg\}
$$
which is a penalized version of the ELBO. Note that taking $(\pi_K)$ as uniform on a finite set $\{1,2,\dots,K_{\max}\}$ leads to the procedure described in \cite{Review}. We discuss below the choice $\pi_K = 2^{-K}$.

\vspace{0.2cm}
We can now state the following result which provides an oracle-type inequality for $\tilde{\pi}_{n,\alpha}^{\hat{K}}(.|X_1^n)$:

\begin{thm}

\label{thm-model-selection}
For any $\alpha \in (0,1)$,
\begin{equation*}
\mathbb{E} \bigg[ \int D_{\alpha}( P_{\theta}, P^0 ) \tilde{\pi}_{n,\alpha}^{\hat{K}}(d\theta|X_1^n) \bigg] \leq \inf_{K \geq 1} \bigg\{ \frac{\alpha}{1-\alpha} \inf_{\theta^* \in \Theta_K(r_{n,K})} \mathcal{K}(P^0,P_{\theta^*}) + \frac{1+\alpha}{1-\alpha}2Kr_{n,K} +  \frac{\log(\frac{1}{\pi_K})}{n(1-\alpha)} \bigg\}.
\end{equation*}
\end{thm}

\vspace{0.2cm}
This oracle inequality shows that our variational distribution adaptively satisfies the best possible balance between bias (misspecification error) and variance (estimation error). If we assume that there is actually a $K_0$ and $\theta^*\in\Theta_{K_0}$ such that $P^0=P_{\theta^*}$ then the theorem will imply
\begin{equation*}
\mathbb{E} \bigg[ \int D_{\alpha}( P_{\theta}, P^0 ) \tilde{\pi}_{n,\alpha}^{\hat{K}}(d\theta|X_1^n) \bigg] \leq \frac{1+\alpha}{1-\alpha}2K_0 r_{n,K_0} +  \frac{\log(\frac{1}{\pi_{K_0}})}{n(1-\alpha)}.
\end{equation*}
Note that this does not mean that $\hat{K}=K_0$, but this means that the convergence rate of $P_\theta$ to $P^0$ $\tilde{\pi}_{n,\alpha}^{\hat{K}}(.|X_1^n)$ is as good as if we actually knew $P_0$. The objective of estimating $K_0$ is a completely different task \cite{yang2005can}. Estimating $K_0$ would also require identifiability conditions that are not necessary for our results.

\vspace{0.2cm}
The variance term is composed of two parts. The first one, $Kr_{n,K}$ up to a multiplicative constant, corresponds to the rate obtained when approximating the true distribution with mixtures of model $\mathcal{M}_K$. The second part of the overall rate can be interpreted as a complexity term over the different models reflecting our prior belief. For instance, if we want to penalize more complex models, we can take $\pi_K=2^{-K}$ and the corresponding term will be of order $K/n$. In practice, as soon as $\frac{1}{n}\lesssim r_{n,K}$, then this penalty term is negligible when compared to the approximating rate $Kr_{n,K}$: this means that this choice can be considered safe, as it does not interfere with the estimation rate.

\vspace{0.2cm}
\section{Conclusion}
\label{section:conclusion}
\vspace{0.2cm}

Using variational inference, we studied consistency of variational approximations for estimation and model selection in mixtures. When considering tempered posteriors, we showed that Variational Bayes is consistent and we gave statistical guarantees to model selection based on the ELBO. For further investigation, it would be interesting to explore the case of Bayesian posteriors when $\alpha=1$. The recent work of Zhang and Gao \cite{Chicago} gives the tools for tackling such an issue, and allows one to consider risk functions different from $\alpha$-Renyi divergence. But the conditions would be more stringent, and misspecification would be more problematic in this case.

Another point of interest is the study of the non-convex optimization program \eqref{defpasELBO}. Indeed, the proposed coordinate optimization can lead to a local extremum, and this implies that one needs to pay attention to initialization. The same problem also occurs in the Expectation-Maximization (EM) algorithm. In practice, users often run EM or CAVI several times with different initial distributions. Many practical ideas were proposed to target the global extrema more efficiently with EM \cite{o2012computational} and could be extended to CAVI. But the question of convergence remains open in theory.

Finally, note that our results are remarkable as there are almost no conditions on the mixtures considered. In this paper we have focused on estimating the true probability distribution $P^0$, even in the well-specified case. We have no results on the estimation of the parameters. In the case of mixtures, these results are extremely difficult to obtain even for Gaussian mixtures \cite{wu2018mixtures}. They require restrictions on the parameters set and lead to different rates of convergence. The consistency of VB for the estimation of the parameters remains open.

\section*{Acknowledgements}

We thank the Associate Editor and the anonymous Referee for their insightful comments on the paper.

\bibliographystyle{plain}

\vspace{0.2cm}
\section{Proofs}
\label{sec:proofs}
\vspace{0.2cm}

\vspace{0.2cm}
\subsection{Some useful lemmas}
\vspace{0.2cm}

We provide in this section two useful lemmas required in many proofs below.

\vspace{0.2cm}
\subsubsection{An upper bound on the Kullback-Leibler divergence between two mixtures}
\vspace{0.2cm}

The lemma below was first stated by \cite{DoWarmuthAndSinger} for mixtures of Gaussians, \cite{DoMinh} checked that the proof remains valid for general mixtures. It is a tool widely used in signal processing \cite{DoLemmaSignalProcessing}. We provide the proof for the sake of completeness.

\begin{lemma}

\label{Do}
Let $p,p^0 \in \mathcal{S}_K$ and $\theta_j,\theta_j^0 \in \Theta$ for $j=1,...,K$. Then,

$$
\mathcal{K}\left( \sum_{j=1}^K p_j^0 q_{\theta_j^0} ,\sum_{j=1}^K p_j q_{\theta_j} \right) \leq \mathcal{K}(p^0,p) + \sum_{j=1}^K p_j^0 \mathcal{K}(q_{\theta_j^0},q_{\theta_j})
$$

\end{lemma}

\vspace{0.2cm}
\begin{proof}
For any nonnegative numbers $\alpha_1,...,\alpha_K$ and positive $\beta_1,...,\beta_K$, we have:

\vspace{0.2cm}
$\begin{aligned}[t]
   \left( \sum_{j=1}^K \alpha_j \right) \log \left( \frac{\sum_{j=1}^K \alpha_j }{\sum_{j=1}^K \beta_j } \right) & = \left( \sum_{j=1}^K \beta_j \right) \left( \frac{\sum_{j=1}^K \alpha_j }{\sum_{j=1}^K \beta_j } \right) \log \left( \frac{\sum_{j=1}^K \alpha_j }{\sum_{j=1}^K \beta_j } \right) \\
                                     & = \left( \sum_{j=1}^K \beta_j \right) \left( \sum_{j=1}^K \frac{\beta_j}{\sum_{l=1}^K \beta_l} \frac{\alpha_j}{\beta_j} \right) \log \left( \sum_{j=1}^K \frac{\beta_j}{\sum_{l=1}^K \beta_l} \frac{\alpha_j}{\beta_j} \right) \\
                                     & = \left( \sum_{j=1}^K \beta_j \right) f\left( \sum_{j=1}^K \frac{\beta_j}{\sum_{l=1}^K \beta_l} \frac{\alpha_j}{\beta_j} \right)
\end{aligned}$

\vspace{0.2cm}
\noindent
where $f$ is the convex function $x \longmapsto x\log(x)$. As $\sum_{j=1}^K \frac{\beta_j}{\sum_{l=1}^K \beta_l}=1$, then using Jensen's inequality:

\vspace{0.2cm}
$\begin{aligned}[t]
   \left( \sum_{j=1}^K \alpha_j \right) \log \left( \frac{\sum_{j=1}^K \alpha_j }{\sum_{j=1}^K \beta_j } \right) & = \left( \sum_{j=1}^K \beta_j \right) f\left( \sum_{j=1}^K \frac{\beta_j}{\sum_{l=1}^K \beta_l} \frac{\alpha_j}{\beta_j} \right) \\
                                     & \leq \left( \sum_{j=1}^K \beta_j \right) \sum_{j=1}^K \frac{\beta_j}{\sum_{l=1}^K \beta_l} f\left( \frac{\alpha_j}{\beta_j}\right) \\
                                     & = \left( \sum_{j=1}^K \beta_j \right) \sum_{j=1}^K \frac{\beta_j}{\sum_{l=1}^K \beta_l}  \frac{\alpha_j}{\beta_j}  \log \left( \frac{\alpha_j}{\beta_j}\right) \\
                                     & = \sum_{j=1}^K \alpha_j\log \left( \frac{\alpha_j}{\beta_j}\right).
\end{aligned}$

\vspace{0.2cm}
\noindent
The inequality remains valid when some or all $\beta_j$'s are zero. Indeed, assume that $\beta_j=0$. If $\alpha_j \ne 0$, then the $j^{th}$ term of the sum in the right-hand side is $\alpha_j \log ( {\alpha_j}/{\beta_j})=+\infty$, and the result is obvious. Otherwise, $\alpha_j=0$, hence the $j^{th}$ term of each sum in the inequality is zero as $\alpha_j \log ( {\alpha_j}/{\beta_j})=0$, and the inequality can be obtained considering only the other numbers.

\vspace{0.2cm}
\noindent
Thus, for $p,p^0 \in \mathcal{S}_K$ and $\theta_j,\theta_j^0 \in \Theta$ for $j=1,...,K$:

\vspace{0.2cm}
$\begin{aligned}[t]
   \mathcal{K}\left( \sum_{j=1}^K p_j^0 q_{\theta_j^0} ,\sum_{j=1}^K p_j q_{\theta_j} \right) & = \int \left( \sum_{j=1}^K p_j^0 q_{\theta_j^0} \right) \log \left( \frac{\sum_{j=1}^K p_j^0 q_{\theta_j^0}}{\sum_{j=1}^K p_j q_{\theta_j}} \right) \\
                                     & \leq \int  \sum_{j=1}^K p_j^0 q_{\theta_j^0}  \log \left( \frac{p_j^0 q_{\theta_j^0}}{p_j q_{\theta_j}} \right) \\
                                     & = \int  \sum_{j=1}^K p_j^0 q_{\theta_j^0}  \log \left( \frac{p_j^0}{p_j} \right) + \int  \sum_{j=1}^K p_j^0 q_{\theta_j^0}  \log \left( \frac{q_{\theta_j^0}}{q_{\theta_j}} \right) \\
                                     & = \sum_{j=1}^K p_j^0 \log \left( \frac{p_j^0}{p_j} \right) \left(\int q_{\theta_j^0} \right) + \sum_{j=1}^K p_j^0 \int  q_{\theta_j^0} \log \left( \frac{q_{\theta_j^0}}{q_{\theta_j}} \right) \\
                                     & = \mathcal{K}(p^0,p) + \sum_{j=1}^K p_j^0 \mathcal{K}(q_{\theta_j^0},q_{\theta_j}),
\end{aligned}$

\vspace{0.2cm}
\noindent which ends the proof.
\end{proof}

\vspace{0.2cm}
\subsubsection{KL-divergence between Gaussian distributions and between Normal-Inverse-Gamma distributions}
\vspace{0.2cm}

We give in this section the Kullback-Leibler divergence between 1-dimensional Gaussian distributions and between Normal-Inverse-Gamma distributions.

\begin{lemma}

\label{thm-kl-divergence-gaussians}
We denote $u$ and $v$ the density functions of the respective Gaussian distributions $\mathcal{N}(\mu_u,\sigma^2_u)$ and $\mathcal{N}(\mu_v,\sigma^2_v)$. Similarly, we denote $p$ and $q$ the two densities of $\mathcal{NIG}(\mu_1,\theta_1^2,a_1,b_1)$ and $\mathcal{NIG}(\mu_2,\theta_2^2,a_2,b_2)$. Then:

$$
\mathcal{K}(u,v) = \frac{1}{2} \log\bigg(\frac{\sigma^2_v}{\sigma_u^2}\bigg) + \frac{\sigma_u^2}{2\sigma_v^2} + \frac{(\mu_v-\mu_u)^2}{2\sigma_v^2} - \frac{1}{2} 
$$
and

\vspace{0.2cm}
$\begin{aligned}[t]
   \mathcal{K}(p,q) = \frac{1}{2} \log\bigg(\frac{\theta_1^2}{\theta_2^2}\bigg) & + \frac{\theta_2^2}{2\theta_1^2} + \frac{\theta_2^2(\mu_2-\mu_1)^2}{2} \frac{a_1}{b_1} - \frac{1}{2}  \\
                                     & + (a_1-a_2) \psi(a_1) + \log\bigg(\frac{\Gamma(a_2)}{\Gamma(a_1)}\bigg) + a_2 \log\bigg(\frac{b_1}{b_2}\bigg) + a_1 \frac{b_2-b_1}{b_1}.
\end{aligned}$
\end{lemma}

\vspace{0.2cm}
\begin{proof}
The first equality is extremely classical so we don't provide the proof. For the second one,

\vspace{0.2cm}
$\begin{aligned}[t]
   \mathcal{K}(p,q) & = \int_{\mathbb{R}_+^*} \int_\mathbb{R} p(x,y) \log\bigg(\frac{p(x,y)}{q(x,y)}\bigg) dx dy  \\
                                     & = \int_{\mathbb{R}_+^*} \int_\mathbb{R} p(x|Y=y) p_Y(y) \log\bigg(\frac{p(x|Y=y)}{q(x|Y=y)} \frac{p_Y(y)}{q_Y(y)}\bigg) dx dy \\
                                     & = \int_{\mathbb{R}_+^*} p_Y(y) \bigg( \int_\mathbb{R} p(x|Y=y) \log\bigg(\frac{p(x|Y=y)}{q(x|Y=y)}\bigg) dx \bigg) dy + \int_{\mathbb{R}_+^*} p_Y(y) \log\bigg(\frac{p_Y(y)}{q_Y(y)}\bigg)dy \\
                                     & = \mathbb{E}_{Y\sim \mathcal{IG}(a_1,b_1)}\bigg[ \mathcal{K}(p(.|Y),q(.|Y)) \bigg] + \mathcal{K}(p_Y,q_Y).
\end{aligned}$

\vspace{0.2cm}
Using the KL-divergence between Gaussians: 
$$\mathcal{K}(p(.|Y),q(.|Y))  = \frac{1}{2} \log\bigg(\frac{\theta_1^2}{\theta_2^2}\bigg) + \frac{\theta_2^2}{2\theta_1^2} + \frac{\theta_2^2(\mu_2-\mu_1)^2}{2Y} - \frac{1}{2}$$
hence
$$
\mathbb{E}_{Y\sim \mathcal{IG}(a_1,b_1)}\bigg[ \mathcal{K}(p(.|Y),q(.|Y))  \bigg] = \frac{1}{2} \log\bigg(\frac{\theta_1^2}{\theta_2^2}\bigg) + \frac{\theta_2^2}{2\theta_1^2} + \frac{\theta_2^2(\mu_2-\mu_1)^2}{2}\mathbb{E}_{Y\sim \mathcal{IG}(a_1,b_1)}\bigg[\frac{1}{Y}\bigg] - \frac{1}{2}
$$
i.e.
$$
\mathbb{E}_{Y\sim \mathcal{IG}(a_1,b_1)}\bigg[ \mathcal{K}(p(.|Y),q(.|Y))  \bigg] = \frac{1}{2} \log\bigg(\frac{\theta_1^2}{\theta_2^2}\bigg) + \frac{\theta_2^2}{2\theta_1^2} + \frac{\theta_2^2(\mu_2-\mu_1)^2}{2} \frac{a_1}{b_1} - \frac{1}{2},
$$
and using the KL-divergence between Inverse-Gamma distributions 
$$\mathcal{K}(p_Y,q_Y) = (a_1-a_2) \psi(a_1) + \log\bigg(\frac{\Gamma(a_2)}{\Gamma(a_1)}\bigg) + a_2 \log\bigg(\frac{b_1}{b_2}\bigg) + a_1 \frac{b_2-b_1}{b_1}$$
where $\Gamma$ and $\psi$ are respectively the Gamma and Digamma functions, we have:

\vspace{0.2cm}
$\begin{aligned}[t]
   \mathcal{K}(p,q) = \frac{1}{2} \log\bigg(\frac{\theta_1^2}{\theta_2^2}\bigg) & + \frac{\theta_2^2}{2\theta_1^2} + \frac{\theta_2^2(\mu_2-\mu_1)^2}{2} \frac{a_1}{b_1} - \frac{1}{2}  \\
                                     & + (a_1-a_2) \psi(a_1) + \log\bigg(\frac{\Gamma(a_2)}{\Gamma(a_1)}\bigg) + a_2 \log\bigg(\frac{b_1}{b_2}\bigg) + a_1 \frac{b_2-b_1}{b_1}.
\end{aligned}$

\end{proof}

\vspace{0.2cm}
\subsection{Proof of Theorem~\ref{thm-general}}
\vspace{0.2cm}

This result relies on an application of Theorem 2.6 in \cite{Tempered} to mixture models. The proof of Theorem 2.6 in \cite{Tempered} itself relies mostly on a deviation inequality from \cite{bhattacharya2016bayesian} and on PAC-Bayesian theory \cite{catoni2004statistical,MR2483528}.

\begin{proof}
Fix $0<\alpha<1$. Theorem 2.6 from \cite{Tempered} gives:
   \begin{align*}
 & \mathbb{E} \bigg[ \int D_{\alpha}\bigg( \sum_{j=1}^K p_j q_{\theta_j} ,\sum_{j=1}^K p_j^0 q_{\theta_j^0} \bigg) \tilde{\pi}_{n,\alpha}(d\theta|X_1^n) \bigg]
 \\
 & \quad \leq \inf_{\rho \in \mathcal{F}} \bigg\{ \frac{\alpha}{1-\alpha} \int \mathcal{K}\bigg(\sum_{j=1}^K p_j^0 q_{\theta_j^0}, \sum_{j=1}^K p_j q_{\theta_j} \bigg) \rho({\rm d} \theta)
 + \frac{\mathcal{K}(\rho,\pi)}{n(1-\alpha)}  \bigg\}.
 \end{align*}
 Thanks to Lemma \ref{Do} 
 $$
 \mathcal{K}\bigg(\sum_{j=1}^K p_j^0 q_{\theta_j^0}, \sum_{j=1}^K p_j q_{\theta_j} \bigg)
 \leq \mathcal{K}(p^0,p) + \sum_{j=1}^K p_j^0 \mathcal{K}(q_{\theta_j^0},q_{\theta_j}).
 $$
 Then
 $$ \mathcal{K}(\rho,\pi) = \mathcal{K}\left(\rho_p \bigotimes_{j=1}^{K}\rho_j,\pi_p \bigotimes_{j=1}^{K}\pi_j\right) = \mathcal{K}(\rho_p,\pi_p) + \sum_{j=1}^K \mathcal{K}(\rho_j,\pi_j) $$
 the last inequality being obtained thanks to Theorem 28 in \cite{van2014renyi}.
Gathering all the pieces together leads to
  \begin{align*}
 & \mathbb{E} \bigg[ \int D_{\alpha}\bigg( \sum_{j=1}^K p_j q_{\theta_j} ,\sum_{j=1}^K p_j^0 q_{\theta_j^0} \bigg) \tilde{\pi}_{n,\alpha}(d\theta|X_1^n) \bigg]
 \\
 & \quad \leq \inf_{\rho \in \mathcal{F}} \bigg\{ \frac{\alpha}{1-\alpha} \bigg[ \int \mathcal{K}(p^0,p) \rho_p(dp) + \sum_{j=1}^K \int \mathcal{K}(q_{\theta_j^0},q_{\theta_j}) \rho_j(d\theta_j) \bigg] + \frac{\mathcal{K}(\rho_p,\pi_p) + \sum_{j=1}^K \mathcal{K}(\rho_j,\pi_j)}{n(1-\alpha)}  \bigg\}
 \end{align*}
 that is the result stated in Theorem~\ref{thm-general}.
\end{proof}

\vspace{0.2cm}
\subsection{Proof of Lemma~\ref{lemma-Dirichlet}}
\vspace{0.2cm}

\begin{proof}
Let us define $\rho_{p,n} \in \mathcal{M}_1^+(\mathcal{S}_K)$ by the following formula $\rho_{p,n}(dp) \propto \mathbf{1}(p \in \mathcal{B}) \pi_p(dp)$ with $$
\mathcal{B} = \bigg\{p \in \mathcal{S}_K / \mathcal{K}(p^0,p) \leq K r_{n,K}'\bigg\}
$$ 
and 
$$
r_{n,K}'=\max\bigg(\frac{1}{K(n-1)},\frac{\log(n(K-1)\Gamma(A)^\frac{K}{K-1}/M_p^0)}{n}\bigg)
$$ where $A=\frac{2}{K}$
and $M_p^0=\max\{p_j^0/j=1,...,K\}$. We adopt the notation $S=\sum_{j=1}^K\alpha_j$ in the following. Recall that by assumption $K \geq 2$ and hence $A=\frac{2}{K}\leq 1$.

\vspace{0.2cm}
First, $\int \mathcal{K}(p^0,p) \rho_{p,n}(dp) \leq K r_{n,K}'$.

\vspace{0.2cm}
Then, let us show that $\mathcal{K}(\rho_{p,n},\pi_p) \leq K n r_{n,K}'$.
For that, let us define 
$$
\mathcal{A} = \bigg\{p \in \mathbb{R}^K / p_j^0 e^{-Kr_{n,K}'} \leq p_j \leq p_j^0 e^{-Kr_{n,K}'} + \frac{p_K^0}{n(K-1)}\hspace{0.2cm} \textnormal{for} \hspace{0.2cm} j=1,...,K-1, \hspace{0.1cm} p_K = 1 - \sum_{j=1}^{K-1} p_j \bigg\}
$$ 
where $K$ is such that $p_K^0 = \max\{p_j^0/j=1,...,K\}$ (this assumption can always be fulfilled by reordering and relabelling the vector components). Then, $p_K^0 \geq \frac{1}{K}$ (otherwise, the sum of the components of $p^0$ would be strictly lower than $1$ and the vector would not be included in $\mathcal{S}_K$). We will show that $\mathcal{A} \subset \mathcal{B}$ and that $\pi_p(\mathcal{A}) \geq e^{-Knr_{n,K}'}$. Then, we will conclude thanks to the following formula: $\mathcal{K}(\rho_{p,n},\pi_p) = -\log(\pi_p(\mathcal{B}))$.

\vspace{0.2cm}
First, let us show that $\mathcal{A} \subset \mathcal{B}$.

\vspace{0.2cm}
Let $p \in \mathcal{A}$. As $p_K = 1 - \sum_{j=1}^{K-1} p_j$, we just need to check that $\mathcal{K}(p^0,p) \leq K r_{n,K}'$ and that $p_j \geq 0$ for each $j=1,...,K$.

\vspace{0.2cm}
The first part can be proven using the definition of $\mathcal{A}$. According to the $K-1$ left-hand side inequalities in the definition of $\mathcal{A}$,

\vspace{0.2cm}
$\begin{aligned}[t]
   \mathcal{K}(p^0,p) = \sum_{j=1}^{K-1} p_j^0 \log\bigg(\frac{p_j^0}{p_j}\bigg) + p_K^0 \log\bigg(\frac{p_K^0}{p_K}\bigg) & \leq \sum_{j=1}^{K-1} p_j^0 \log(e^{Kr_{n,K}'}) + p_K^0 \log\bigg(\frac{p_K^0}{p_K}\bigg) \\
                                     & = \sum_{j=1}^{K-1} p_j^0 Kr_{n,K}' + p_K^0 \log\left(\frac{p_K^0}{p_K}\right) \\
                       & = (1-p_K^0) Kr_{n,K}' + p_K^0 \log\left(\frac{p_K^0}{p_K}\right).
\end{aligned}$

\vspace{0.2cm}
All we need to show now is that $\log\left(\frac{p_K^0}{p_K}\right) \leq Kr_{n,K}'$. This comes from the following inequalities:

\vspace{0.2cm}
$\begin{aligned}[t]
   \log\bigg(\frac{p_K^0}{p_K}\bigg) = \log\bigg(\frac{p_K^0}{1 - \sum_{j=1}^{K-1} p_j}\bigg) & \leq \log\bigg(\frac{p_K^0}{1 - \sum_{j=1}^{K-1} p_j^0 e^{-Kr_{n,K}'} - \frac{p_K^0}{n}} \bigg)  \\
                                     & = \log\bigg(\frac{p_K^0}{1 - (1-p_K^0) e^{-Kr_{n,K}'} - \frac{p_K^0}{n}} \bigg)  \\
                                     & \leq \frac{p_K^0}{1 - (1-p_K^0) e^{-Kr_{n,K}'} - \frac{p_K^0}{n}} - 1 \\
                                     & = \frac{p_K^0 - 1 + (1-p_K^0) e^{-Kr_{n,K}'} + \frac{p_K^0}{n} }{1 - (1-p_K^0) e^{-Kr_{n,K}'} - \frac{p_K^0}{n}}
\end{aligned}$

i.e.
\vspace{0.2cm}
$\begin{aligned}[t]
   \log\bigg(\frac{p_K^0}{p_K}\bigg) & \leq \frac{p_K^0 - 1 + (1-p_K^0) e^{-Kr_{n,K}'} + \frac{p_K^0}{n} }{1 - (1-p_K^0) e^{-Kr_{n,K}'} - \frac{p_K^0}{n}} && \\                      & = \frac{\frac{p_K^0}{n} - (1-p_K^0)(1-e^{-Kr_{n,K}'}) }{p_K^0 (1-\frac{1}{n}) + (1-p_K^0)(1-e^{-Kr_{n,K}'})} && \\                                     & = \frac{\frac{1}{n} - (\frac{1}{p_K^0}-1)(1-e^{-Kr_{n,K}'}) }{(1-\frac{1}{n}) + (\frac{1}{p_K^0}-1)(1-e^{-Kr_{n,K}'})} && \\
                                   & \leq \frac{\frac{1}{n}}{1-\frac{1}{n}} = \frac{1}{n-1} && \\
                                   & \leq Kr_{n,K}'. &&
\end{aligned}$

\vspace{0.2cm}
Hence $ \hspace{0.2cm}
\mathcal{K}(p^0,p) \leq  (1-p_K^0) Kr_{n,K}' + p_K^0 \log\left(\frac{p_K^0}{p_K}\right) \leq  (1-p_K^0) Kr_{n,K}' + p_K^0 Kr_{n,K}' = Kr_{n,K}'$.

\vspace{0.2cm}
On the other hand, for $j=1,...,K-1$, $p_j \geq p_j^0 e^{-Kr_{n,K}'} \geq 0$ and:

$\begin{aligned}[t]
   p_K = 1 - \sum_{j=1}^{K-1} p_j & \geq 1 - \sum_{j=1}^{K-1} \big( p_j^0 e^{-Kr_{n,K}'} + \frac{p^0_K}{n(K-1)} \big)  \\
                                     & = 1 - \bigg( ( 1 - p_K^0 ) e^{-Kr_{n,K}'} + \frac{p_K^0}{n} \bigg) \\
                                     &  \geq 1 - ( 1 - p_K^0 ) e^{-Kr_{n,K}'} - p_K^0 \\
                                     & = ( 1 - p_K^0 )( 1 - e^{-Kr_{n,K}'})\\
                                     & \geq 0.
\end{aligned}$

\vspace{0.2cm}
Then, $p \in \mathcal{B}$, and finally $\mathcal{A} \subset \mathcal{B}$.

\vspace{0.2cm}
Now, let us show that $\pi_p(\mathcal{A}) \geq e^{-Knr_{n,K}'}$.

\vspace{0.2cm}
Let us denote $f$ the density of the $\pi_p=\mathcal{D}_K(\alpha_1,...,\alpha_K)$ Dirichlet distribution:
\[
f \left(p\right) = \frac{\Gamma\big(S\big)}{\prod\limits_{j=1}^K \Gamma(\alpha_j)} \prod_{j=1}^K p_j^{\alpha_j - 1} \hspace{0.1cm} \mathbf{1}(p \in \mathcal{S}_K).
\]

\vspace{0.2cm}
Thus, we can lower bound $\pi_p(\mathcal{A})$:

$\begin{aligned}[t]
   \pi_p(\mathcal{A}) = \int_{\mathcal{A}} f(p_1,...,p_K) dp & = \int_{\mathcal{A}}  \frac{\Gamma\big(S\big)}{\prod\limits_{j=1}^K \Gamma(\alpha_j)} \prod_{j=1}^K p_j^{\alpha_j - 1} \hspace{0.1cm} \mathbf{1}(p \in \mathcal{S}_K) dp  \\
                                     & \geq \frac{\Gamma\big(S\big)}{\prod\limits_{j=1}^K \Gamma(\alpha_j)} \prod_{j=1}^{K-1} \int_{p_j^0 e^{-Kr_{n,K}'}}^{p_j^0 e^{-Kr_{n,K}'} + \frac{p_K^0}{n(K-1)} }   p_j^{\alpha_j - 1} \hspace{0.1cm} dp_j 
\end{aligned}$

\noindent as for $p \in \mathcal{A}$, $0 \leq p_j^0 e^{-Kr_{n,K}'} \leq p_j \leq p_j^0 e^{-Kr_{n,K}'} + \frac{p_K^0}{n(K-1)} \leq 1 $ for each $j=1,...,K-1$ (as $\mathcal{A} \subset \mathcal{B}$), and then $p_j^{\alpha_j-1} \geq 1$.

\vspace{0.2cm}
Then, by definition of $r_{n,K}'$, $\frac{p_K^0}{n(K-1)} \geq \Gamma(A)^{\frac{K}{K-1}} e^{-nr_{n,K}'}$, and using inequalities $\Gamma(A)\geq\Gamma(\alpha_j)$ as $A \leq \alpha_j \leq 1$ and $\Gamma(S)\geq1$ as $S\geq2$,

\vspace{0.2cm}
$\begin{aligned}[t]
   \pi_p(\mathcal{A}) & \geq \frac{\Gamma\big(S\big)}{\prod\limits_{j=1}^K \Gamma(\alpha_j)} \prod_{j=1}^{K-1} \int_{p_j^0 e^{-Kr_{n,K}'}}^{p_j^0 e^{-Kr_{n,K}'} + \Gamma(A)^{\frac{K}{K-1}} e^{-nr_{n,K}'}} p_j^{\alpha_j - 1} \hspace{0.1cm} dp_j \\
                                     & \geq \frac{\Gamma\big(S\big)}{\prod\limits_{j=1}^K \Gamma(\alpha_j)} \prod_{j=1}^{K-1} \int_{p_j^0 e^{-Kr_{n,K}'}}^{p_j^0 e^{-Kr_{n,K}'} + \Gamma(A)^{\frac{K}{K-1}} e^{-nr_{n,K}'} } \hspace{0.1cm} dp_j \\
                                     & = \frac{\Gamma\big(S\big)}{\prod\limits_{j=1}^K \Gamma(\alpha_j)} \prod_{j=1}^{K-1}  \Gamma(A)^{\frac{K}{K-1}} e^{-nr_{n,K}'} \\
                                     & = \frac{\Gamma\big(S\big)}{\prod\limits_{j=1}^K \Gamma(\alpha_j)}  \Gamma(A)^K e^{-n(K-1)r_{n,K}'} \\
                                    & \geq
                                    e^{-nKr_{n,K}'}.
\end{aligned}$

\vspace{0.2cm}
Hence, as $\mathcal{A} \subset \mathcal{B}$, $\pi_p(\mathcal{B}) \geq \pi_p(\mathcal{A}) \geq e^{-nKr_{n,K}'}$, and finally, $\mathcal{K}(\rho_{p,n},\pi_p) = -\log(\pi_p(\mathcal{B})) \leq Knr_{n,K}'$.

\vspace{0.2cm}
We just proved the lemma but with the rate $r_{n,K}'$ instead of the value $r_{n,K}$ used in the lemma. We can conlude by noticing that the result is valid for every $r$ such that $r_{n,K}' \leq r$, and that in particular $r_{n,K}' \leq r_{n,K}$. This last result comes from the inequality:
$$ \Gamma(A) \leq \frac{\Gamma(1+\frac{A}{2})}{\left(\frac{A}{2}\right)^{1-\frac{A}{2}}}$$
which is a direct application of the left-hand side of inequality (3.2) part 3 in \cite{Gamma} with $x=\frac{A}{2}>0$ and $\lambda=\frac{A}{2} \in (0,1)$.
As, $1+\frac{A}{2} \in [1,2]$, then $\Gamma(1+\frac{A}{2}) \leq 1$, and $\frac{1}{\left(\frac{A}{2}\right)^{1-\frac{A}{2}}}=K^{1-\frac{A}{2}} \leq K$. Thus: 
$$ \Gamma(A) \leq K $$
and as $K \geq 2$ and $p^0_K\geq\frac{1}{K}$, it follows that 
$$ \log((K-1)\Gamma(A)^\frac{K}{K-1}/p_K^0) \leq \log(K(K)^\frac{K}{K-1}K) \leq \log(K(K)^2K) \leq \log(K^4) $$
i.e. $r_{n,K}' \leq \max(\frac{1}{K(n-1)},\frac{\log(nK^4)}{n}) \leq \max(\frac{1}{K(n-1)},\frac{4\log(nK)}{n}) $. Besides, $\frac{n}{n-1}=1+\frac{1}{n-1}\leq 2$ implies $\frac{1}{K(n-1)} \leq \frac{1}{2(n-1)}\leq \frac{1}{n} \leq \frac{4\log(2)}{n} \leq \frac{4\log(nK)}{n}$, and finally $r_{n,K}' \leq \frac{4\log(nK)}{n}=r_{n,K} $.

\end{proof}

\vspace{0.2cm}
\subsection{Proof of Corollary~\ref{thm-multinomial}}
\vspace{0.2cm}

\begin{proof}

According to Lemma \ref{lemma-Dirichlet}, there exists a distribution $\rho_{p,n} \in \mathcal{M}_1^+(\mathcal{S}_K)$ such that
\begin{equation*}
  \int \mathcal{K}(p^0,p) \rho_{p,n}(dp) \leq K  \frac{4\log(nK)}{n}
\end{equation*}
and 
\begin{equation*}
  \mathcal{K}(\rho_{p,n},\pi_p) \leq K n  \frac{4\log(nK)}{n}.
\end{equation*}

\vspace{0.2cm}
Similarly, the same result states that there exists distributions $\rho_{j,n} \in \mathcal{M}_1^+(\mathcal{S}_V)$ for $j=1,...,K$ such that
\begin{equation*}
  \int \mathcal{K}(q_{\theta_j^0},q_{\theta_j}) \rho_{j,n}(d\theta_j) \leq \frac{4V\log(nV)}{n}
\end{equation*}
and 
\begin{equation*}
  \mathcal{K}(\rho_{j,n},\pi_j) \leq n \frac{4V\log(nV)}{n}.
\end{equation*}

\vspace{0.2cm}
We conclude using Theorem \ref{thm-general}:
\begin{equation*}
   \begin{split}
       \mathbb{E} \bigg[ \int D_{\alpha}( \sum_{j=1}^K p_j q_{\theta_j} ,\sum_{j=1}^K & p_j^0 q_{\theta_j^0} ) \tilde{\pi}_{n,\alpha}(d\theta|X_1^n) \bigg] \leq \frac{1+\alpha}{1-\alpha} \bigg[ \frac{8KV\log(nV)}{n} \bigvee \frac{8K\log(nK)}{n} \bigg].
   \end{split} 
\end{equation*}

\end{proof}

\vspace{0.2cm}
\subsection{Proof of Corollary~\ref{thm-Gaussian}}
\vspace{0.2cm}

\begin{proof}
For $R_{j,n}=\frac{1}{n} \bigvee \frac{1}{n} \bigg[ \frac{1}{2} \log\bigg(\frac{n}{2}\bigg) + \frac{V^2}{n\mathcal{V}^2} + \log\bigg(\frac{\mathcal{V}}{V}\bigg) + \frac{(\mu_j^0)^2}{2\mathcal{V}^2} - \frac{1}{2} \bigg]$ (for $j=1,...,K$), there exists distributions $\rho_{j,n} \in \mathcal{M}_1^+(\mathcal{S}_K)$ for $j=1,...,K$ such that
\begin{equation*}
  \int \mathcal{K}(q_{\mu_j^0},q_{\mu_j}) \rho_{j,n}(d\mu_j) \leq R_{j,n}
\end{equation*}
and 
\begin{equation*}
  \mathcal{K}(\rho_{j,n},\pi_j) \leq n R_{j,n}.
\end{equation*}

\vspace{0.2cm}
Indeed, let us define $\rho_{j,n}$ as a Gaussian distribution of mean $\mu_j^0$ and variance $\frac{2V^2}{n}$. According to Lemma \ref{thm-kl-divergence-gaussians}:
\[
\mathcal{K}(q_{\mu_j^0},q_{\mu_j}) = \frac{(\mu_j-\mu_j^0)^2}{2V^2}.
\]

\vspace{0.2cm}
Then,
$\begin{aligned}[t]
   \int \mathcal{K}(q_{\mu_j^0},q_{\mu_j}) \rho_{j,n}(d\mu_j) & = \frac{1}{2V^2}  \mathbb{E}_{\mu_j\sim\rho_{j,n}}[(\mu_j-\mu_j^0)^2] \\
                                     & = \frac{1}{2V^2} \times \frac{2V^2}{n} \\
                                     & = \frac{1}{n} \\
                                     & \leq R_{j,n}.
\end{aligned}$

\vspace{0.2cm}
We can apply Lemma \ref{thm-kl-divergence-gaussians} again to conclude: 

\vspace{0.2cm}
$\begin{aligned}[t]
   \mathcal{K}(\rho_{j,n},\pi_j) & = \frac{1}{2} \log\bigg(\frac{n\mathcal{V}^2}{2V^2}\bigg) + \frac{V^2}{n\mathcal{V}^2} + \frac{(\mu_j^0)^2}{2\mathcal{V}^2} - \frac{1}{2}  \\
                                     & = \frac{1}{2} \log\bigg(\frac{n}{2}\bigg) + \frac{V^2}{n\mathcal{V}^2} + \log\bigg(\frac{\mathcal{V}}{V}\bigg) + \frac{(\mu_j^0)^2}{2\mathcal{V}^2} - \frac{1}{2} \\
                                     & = n \times \frac{1}{n} \bigg[ \frac{1}{2} \log\bigg(\frac{n}{2}\bigg) + \frac{V^2}{n\mathcal{V}^2} + \log\bigg(\frac{\mathcal{V}}{V}\bigg) + \frac{(\mu_j^0)^2}{2\mathcal{V}^2} - \frac{1}{2} \bigg] \\
                                     & \leq n R_{j,n}.
\end{aligned}$

\vspace{0.2cm}
In addition, Lemma \ref{lemma-Dirichlet} tells us that there exists a distribution $\rho_{p,n} \in \mathcal{M}_1^+(\mathcal{S}_K)$ such that
\begin{equation*}
  \int \mathcal{K}(p^0,p) \rho_{p,n}(dp) \leq K \frac{4\log(nK)}{n} 
\end{equation*}
and 
\begin{equation*}
  \mathcal{K}(\rho_{p,n},\pi_p) \leq  n K  \frac{4\log(nK)}{n}.
\end{equation*}

\vspace{0.2cm}
For $r_{n,K} = \frac{4\log(nK)}{n} \bigvee_{j=1}^{K} R_{j,n} =  \frac{4\log(nK)}{n} \bigvee \frac{1}{n}  \bigvee_{j=1}^{K} \frac{1}{n} \bigg[ \frac{1}{2} \log\bigg(\frac{n}{2}\bigg) + \frac{V^2}{n\mathcal{V}^2} + \log\bigg(\frac{\mathcal{V}}{V}\bigg) + \frac{(\mu_j^0)^2}{2\mathcal{V}^2} - \frac{1}{2} \bigg]$ i.e. $r_{n,K} = \frac{4\log(nK)}{n}  \bigvee_{j=1}^{K} \frac{1}{n} \bigg[ \frac{1}{2} \log\bigg(\frac{n}{2}\bigg) + \frac{V^2}{n\mathcal{V}^2} + \log\bigg(\frac{\mathcal{V}}{V}\bigg) + \frac{(\mu_j^0)^2}{2\mathcal{V}^2} - \frac{1}{2} \bigg]$, we finally obtain the required inequality using Theorem \ref{thm-general}.

\end{proof}

\vspace{0.2cm}
\subsection{Proof of Corollary~\ref{thm-Gaussian-variance}}
\vspace{0.2cm}

\vspace{0.2cm}
\subsubsection{Normal-Inverse-Gamma prior}
\vspace{0.2cm}
\begin{proof}
First, let us focus on the first result, when the chosen prior is the Normal-Inverse-Gamma $\pi_j = \mathcal{NIG}(0,\mathcal{V}^{-2},1,\gamma^2)$ for each $j=1,...,K$. In order to obtain the required rate 
$$
r_{n,K} =  \frac{4\log(nK)}{n} \bigvee_{j=1}^{K} \frac{1}{n} \bigg[  2 \log(n\sqrt{\mathcal{V}}) + \frac{1}{2n\mathcal{V}^2} + \frac{(\mu_j^0)^2}{2(\sigma_j^0)^2\mathcal{V}^2} + \log\bigg(\frac{(\sigma_j^0)^2}{\gamma^2}\bigg) + \frac{\gamma^2}{(\sigma_j^0)^2}  - \frac{1}{2}\log(2\pi)  \bigg],
$$
we proceed as previously and find a variational density on both the mean and the variance such that the two different terms $\int \mathcal{K}(q_{(\mu_j^0,(\sigma_j^0)^2)},q_{(\mu_j,\sigma_j^2)}) \rho_{j,n}(d\mu_j,d\sigma_j^2)$ and $\mathcal{K}(\rho_{j,n},\pi_j)$ are upper bounded for $j=1,...,K$.

\vspace{0.2cm}
Let us define $\rho_{j,n}$ as a Normal-Inverse-Gamma distribution $\mathcal{NIG}(\mu_j^0,\lambda_n,a_n,b_n)$ where $\lambda_n$, $a_n$ and $b_n$ are hyperparameters that we will make precise later. Using Lemma \ref{thm-kl-divergence-gaussians}:
\[
\mathcal{K}(q_{(\mu_j^0,(\sigma_j^0)^2)},q_{(\mu_j,\sigma_j^2)}) = \frac{1}{2} \log\bigg(\frac{\sigma^2_j}{(\sigma_j^0)^2}\bigg) + \frac{(\sigma_j^0)^2}{2\sigma_j^2} + \frac{(\mu_j-\mu_j^0)^2}{2\sigma_j^2} -\frac{1}{2}.
\]

\vspace{0.2cm}
Then, $\begin{aligned}[t]
   \int \mathcal{K}(q_{(\mu_j^0,(\sigma_j^0)^2)},q_{(\mu_j,\sigma_j^2)}) \rho_{j,n}(d\mu_j) = \frac{1}{2} \mathbb{E}_{(\mu_j,\sigma_j^2)\sim\rho_{j,n}}\bigg[&\log\bigg(\frac{\sigma^2_j}{(\sigma_j^0)^2}\bigg)\bigg] + \mathbb{E}_{(\mu_j,\sigma_j^2)\sim\rho_{j,n}}\bigg[\frac{(\sigma_j^0)^2}{2\sigma^2_j}\bigg] \\ 
& + \mathbb{E}_{(\mu_j,\sigma_j^2)\sim\rho_{j,n}}\bigg[\frac{(\mu_j-\mu_j^0)^2}{2\sigma^2_j}\bigg] - \frac{1}{2}.
\end{aligned}$

\vspace{0.2cm}
As $$\mathbb{E}_{(\mu_j,\sigma_j^2)\sim\rho_{j,n}}\bigg[\frac{(\sigma_j^0)^2}{2\sigma^2_j}\bigg] = \frac{(\sigma_j^0)^2}{2} \mathbb{E}_{(\mu_j,\sigma_j^2)\sim\rho_{j,n}}\bigg[\frac{1}{\sigma^2_j}\bigg] = \frac{(\sigma_j^0)^2}{2} \frac{a_n}{b_n},$$

$$\frac{1}{2} \mathbb{E}_{(\mu_j,\sigma_j^2)\sim\rho_{j,n}}\bigg[\log\bigg(\frac{\sigma^2_j}{(\sigma_j^0)^2}\bigg)\bigg] = \frac{1}{2} \big( \log(b_n) - \psi(a_n) \big) - \frac{1}{2}\log((\sigma_j^0)^2)$$

and $$\mathbb{E}_{(\mu_j,\sigma_j^2)\sim\rho_{j,n}}\bigg[\frac{(\mu_j-\mu_j^0)^2}{2\sigma_j^2}\bigg] = \mathbb{E}_{\sigma_j^2\sim \mathcal{IG}(a_n,b_n) }\bigg[ \frac{1}{2\sigma_j^2} . \mathbb{E}_{\mu_j\sim\mathcal{N}(\mu_j^0,\frac{\sigma_j^2}{\lambda_n}) } [ (\mu_j-\mu_j^0)^2 ] \bigg] $$
i.e.
$$\mathbb{E}_{(\mu_j,\sigma_j^2)\sim\rho_{j,n}}\bigg[\frac{(\mu_j-\mu_j^0)^2}{2\sigma_j^2}\bigg] = \mathbb{E}_{\sigma_j^2\sim \mathcal{IG}(a_n,b_n) }\bigg[ \frac{1}{2\sigma_j^2} . \frac{\sigma_j^2}{\lambda_n} \bigg] = \frac{1}{2\lambda_n},$$
we get: 
$$
\int \mathcal{K}(q_{(\mu_j^0,(\sigma_j^0)^2)},q_{(\mu_j,\sigma_j^2)}) \rho_{j,n}(d\mu_j) = - \frac{1}{2} + \frac{(\sigma_j^0)^2}{2} \frac{a_n}{b_n} + \frac{1}{2\lambda_n} + \frac{1}{2} \big( \log(b_n) - \psi(a_n) \big) - \frac{1}{2}\log((\sigma_j^0)^2).
$$

\vspace{0.2cm}
Now, we compute the term $\mathcal{K}(\rho_{j,n},\pi_j)$ using the fomula giving the Kullback-Leibler divergence between two Gaussian-Inverse-Gamma distributions. Using Lemma \ref{thm-kl-divergence-gaussians}:

\vspace{0.2cm}
$\begin{aligned}[t]
   \mathcal{K}(\rho_{j,n},\pi_j) = \frac{1}{2} \log\bigg(\frac{\lambda_n}{\mathcal{V}^{-2}}\bigg) & + \frac{\mathcal{V}^{-2}}{2\lambda_n} + \frac{\mathcal{V}^{-2}(\mu_j^0)^2}{2} \frac{a_n}{b_n} - \frac{1}{2}  \\
                                     & + (a_n-1) \psi(a_n) + \log\bigg(\frac{1}{\Gamma(a_n)}\bigg) + \log\bigg(\frac{b_n}{\gamma^2}\bigg) + a_n \frac{\gamma^2-b_n}{b_n}.
\end{aligned}$

\vspace{0.2cm}
Then, for $\lambda_n=n$, $a_n=n$ and $b_n=n(\sigma_j^0)^2$:

\[
\int \mathcal{K}(q_{(\mu_j^0,(\sigma_j^0)^2)},q_{(\mu_j,\sigma_j^2)}) \rho_{j,n}(d\mu_j) = \frac{1}{2n} + \frac{1}{2} \big(\log(n) - \psi(n)\big) \leq \frac{1}{2n}+\frac{1}{4n} + \frac{1}{24n^2} \leq R_{j,n}
\]

and

$\begin{aligned}[t]
   \mathcal{K}(\rho_{j,n},\pi_j) & =  \frac{1}{2} \log(n\mathcal{V}^{2}) + \frac{1}{2n\mathcal{V}^2} + \frac{(\mu_j^0)^2}{2(\sigma_j^0)^2\mathcal{V}^2} - \frac{1}{2} + \log\bigg(\frac{(\sigma_j^0)^2}{\gamma^2}\bigg) + \frac{\gamma^2-n(\sigma_j^0)^2}{(\sigma_j^0)^2}
 \\
                                     & \hspace{2cm}  + (n-1)\psi(n) + \log(n) - \log\Gamma(n) \\
                                      & \leq \frac{1}{2} \log(n\mathcal{V}^{2}) + \frac{1}{2n\mathcal{V}^2} + \frac{(\mu_j^0)^2}{2(\sigma_j^0)^2\mathcal{V}^2} - \frac{1}{2} + \log\bigg(\frac{(\sigma_j^0)^2}{\gamma^2}\bigg) + \frac{\gamma^2}{(\sigma_j^0)^2}-n\\
                & \hspace{2cm}  + n\psi(n) + \log(n) - \log(n-1)! \\   
                & \leq \frac{1}{2} \log(n\mathcal{V}^{2}) + \frac{1}{2n\mathcal{V}^2} + \frac{(\mu_j^0)^2}{2(\sigma_j^0)^2\mathcal{V}^2} - \frac{1}{2} + \log\bigg(\frac{(\sigma_j^0)^2}{\gamma^2}\bigg) + \frac{\gamma^2}{(\sigma_j^0)^2}-n\\
                & \hspace{2cm}  + \bigg( n\log(n)-\frac{n}{2n}-\frac{n}{12n^2}+\frac{n}{120n^4} \bigg) + \log(n) \\
                & \hspace{2cm} +\bigg( - \frac{1}{2}\log(2\pi) +n-1 -n\log(n-1) +\frac{1}{2} \log(n-1) \bigg) \\
                & \leq \frac{1}{2} \log(n\mathcal{V}^{2}) + \frac{1}{2n\mathcal{V}^2} + \frac{(\mu_j^0)^2}{2(\sigma_j^0)^2\mathcal{V}^2} - \frac{3}{2} + \log\bigg(\frac{(\sigma_j^0)^2}{\gamma^2}\bigg) + \frac{\gamma^2}{(\sigma_j^0)^2}\\
                & \hspace{2cm}  + n\log\bigg(\frac{n}{n-1}\bigg)-\frac{1}{2} + \log(n) - \frac{1}{2}\log(2\pi) +\frac{1}{2} \log(n-1) \\
                & \leq \frac{1}{2n\mathcal{V}^2} + \frac{(\mu_j^0)^2}{2(\sigma_j^0)^2\mathcal{V}^2} - \frac{3}{2} + \bigg( \log\bigg(\frac{(\sigma_j^0)^2}{\gamma^2}\bigg) + \frac{\gamma^2}{(\sigma_j^0)^2}  - \frac{1}{2}\log(2\pi) \bigg)\\
                & \hspace{2cm}  + 2 -\frac{1}{2} + \bigg( \frac{1}{2} \log(n\mathcal{V}^{2}) + \frac{3}{2} \log(n) \bigg)\\
                & = \frac{1}{2n\mathcal{V}^2} + \frac{(\mu_j^0)^2}{2(\sigma_j^0)^2\mathcal{V}^2} + \bigg( \log\bigg(\frac{(\sigma_j^0)^2}{\gamma^2}\bigg) + \frac{\gamma^2}{(\sigma_j^0)^2}  - \frac{1}{2}\log(2\pi) \bigg) + 2 \log(n\sqrt{\mathcal{V}}) \\
                & = n \times \frac{1}{n} \bigg[  2 \log(n\sqrt{\mathcal{V}}) + \frac{1}{2n\mathcal{V}^2} + \frac{(\mu_j^0)^2}{2(\sigma_j^0)^2\mathcal{V}^2} + \log\bigg(\frac{(\sigma_j^0)^2}{\gamma^2}\bigg) + \frac{\gamma^2}{(\sigma_j^0)^2}  - \frac{1}{2}\log(2\pi) \bigg] \\
                                     & \leq n R_{j,n}
\end{aligned}$

with $R_{j,n}=\frac{1}{n} \bigvee \frac{1}{n} \bigg[  2 \log(n\sqrt{\mathcal{V}}) + \frac{1}{2n\mathcal{V}^2} + \frac{(\mu_j^0)^2}{2(\sigma_j^0)^2\mathcal{V}^2} + \log\bigg(\frac{(\sigma_j^0)^2}{\gamma^2}\bigg) + \frac{\gamma^2}{(\sigma_j^0)^2}  - \frac{1}{2}\log(2\pi) \bigg]$
where we used Theorem 5 in \cite{UpperBoundsDigamma} and inequality (1.15) in \cite{LowerBoundFactorial}:
$$ \forall t>0, \hspace{0.2cm} \log(t)-\frac{1}{2t}-\frac{1}{12t^2} < \psi(t) < \log(t)-\frac{1}{2t}-\frac{1}{12t^2}+\frac{1}{120t^4}, $$
$$ \forall n\geq2, \hspace{0.5cm} n! > \sqrt{2\pi} e^{-n} n^{n+1/2}. $$

\vspace{0.4cm}

Recall again that by Lemma \ref{lemma-Dirichlet}, there exists a distribution $\rho_{p,n} \in \mathcal{M}_1^+(\mathcal{S}_K)$ such that
\begin{equation*}
  \int \mathcal{K}(p^0,p) \rho_{p,n}(dp) \leq K \frac{4\log(nK)}{n}
\end{equation*}
and 
\begin{equation*}
  \mathcal{K}(\rho_{p,n},\pi_p) \leq  n K \frac{4\log(nK)}{n}.
\end{equation*}

We can finally conclude using again Theorem \ref{thm-general} with $ r_{n,K}=\frac{4\log(nK)}{n} \bigvee_{j=1}^{K} R_{j,n} $ i.e.
$$
r_{n,K}=  \frac{4\log(nK)}{n} \bigvee_{j=1}^{K} \frac{1}{n} \bigg[  2 \log(n\sqrt{\mathcal{V}}) + \frac{1}{2n\mathcal{V}^2} + \frac{(\mu_j^0)^2}{2(\sigma_j^0)^2\mathcal{V}^2} + \log\bigg(\frac{(\sigma_j^0)^2}{\gamma^2}\bigg) + \frac{\gamma^2}{(\sigma_j^0)^2}  - \frac{1}{2}\log(2\pi)  \bigg].
$$

\end{proof}

\vspace{0.2cm}
\subsubsection{Factorized prior}
\vspace{0.2cm}

\begin{proof}
Let us focus now on the case of independant priors $\pi_j = \mathcal{N}(0,\mathcal{V}^2) \bigotimes \mathcal{IG}(1,\gamma^2)$ for $j=1,...,K$. The proof is almost the same as previously.

\vspace{0.2cm}
We define here $\rho_{j,n}$ as the product measure of Normal distribution $\mathcal{N}(\mu_j^0,\theta_n^2)$ and of an Inverse-Gamma distribution $\mathcal{IG}(a_n,b_n)$ where $\theta_n^2$, $a_n$ and $b_n$ are hyperparameters to be described later. Then, we have again:
\[
\mathcal{K}(q_{(\mu_j^0,(\sigma_j^0)^2)},q_{(\mu_j,\sigma_j^2)}) = \frac{1}{2} \log\bigg(\frac{\sigma^2_j}{(\sigma_j^0)^2}\bigg) + \frac{(\sigma_j^0)^2}{2\sigma_j^2} + \frac{(\mu_j-\mu_j^0)^2}{2\sigma_j^2} -\frac{1}{2}.
\]

Hence,

\vspace{0.2cm}
$\begin{aligned}[t]
   \int \mathcal{K}(q_{(\mu_j^0,(\sigma_j^0)^2)},q_{(\mu_j,\sigma_j^2)}) \rho_{j,n}(d\mu_j) = \frac{1}{2} \mathbb{E}&_{\sigma_j^2\sim\mathcal{IG}(a_n,b_n)}\bigg[\log\bigg(\frac{\sigma^2_j}{(\sigma_j^0)^2}\bigg)\bigg] + \mathbb{E}_{\sigma_j^2\sim\mathcal{IG}(a_n,b_n)}\bigg[\frac{(\sigma_j^0)^2}{2\sigma_j^2}\bigg] \\ & + \mathbb{E}_{\mu_j\sim\mathcal{N}(\mu_j^0,\theta_n^2)}\big[(\mu_j-\mu_j^0)^2\big] \mathbb{E}_{\sigma_j^2\sim\mathcal{IG}(a_n,b_n)}\bigg[\frac{1}{2\sigma_j^2}\bigg] -\frac{1}{2}.
\end{aligned}$

i.e. $$\int \mathcal{K}(q_{(\mu_j^0,(\sigma_j^0)^2)},q_{(\mu_j,\sigma_j^2)}) \rho_{j,n}(d\mu_j) = -\frac{1}{2} + \frac{a_n}{2b_n} \big( (\sigma_j^0)^2 + \theta_n^2 \big) + \frac{1}{2} \big( \log(b_n) - \psi(a_n) \big) - \frac{1}{2}\log((\sigma_j^0)^2).$$

\vspace{0.2cm}
Then we compute the term $\mathcal{K}(\rho_{j,n},\pi_j)$ as the sum of the Kullback-Leibler divergence between two Gaussian distributions and between two Inverse-Gamma distributions:

\vspace{0.2cm}
$\begin{aligned}[t]
   \mathcal{K}(\rho_{j,n},\pi_j) = \frac{1}{2} \log\bigg(\frac{\mathcal{V}^2}{\theta_n^2}\bigg) & + \frac{\theta_n^2}{2\mathcal{V}^2} + \frac{(\mu_j^0)^2}{2\mathcal{V}^2} - \frac{1}{2}  \\
                                     & + (a_n-1) \psi(a_n) + \log\bigg(\frac{1}{\Gamma(a_n)}\bigg) + \log\bigg(\frac{b_n}{\gamma^2}\bigg) + a_n \frac{\gamma^2-b_n}{b_n}.
\end{aligned}$

\vspace{0.4cm}

Then, for $\theta_n^2=\frac{(\sigma_j^0)^2}{n}$, $a_n=
n$ and $b_n=n(\sigma_j^0)^2$:

\[
\int \mathcal{K}(q_{(\mu_j^0,(\sigma_j^0)^2)},q_{(\mu_j,\sigma_j^2)}) \rho_{j,n}(d\mu_j) = \frac{1}{2n} + \frac{1}{2} \big(\log(n) - \psi(n)\big) \leq \frac{1}{2n}+\frac{1}{4n} + \frac{1}{24n^2} \leq R_{j,n}
\]

and

\vspace{0.2cm}
$\begin{aligned}[t]
   \mathcal{K}(\rho_{j,n},\pi_j) & = \frac{1}{2} \log\bigg(\frac{n\mathcal{V}^2}{(\sigma_j^0)^2}\bigg) + \frac{(\sigma_j^0)^2}{2n\mathcal{V}^2} + \frac{(\mu_j^0)^2}{2\mathcal{V}^2} - \frac{1}{2} + \log\bigg(\frac{(\sigma_j^0)^2}{\gamma^2}\bigg) + \frac{\gamma^2-n(\sigma_j^0)^2}{(\sigma_j^0)^2}
 \\
                                     & \hspace{2cm}  + (n-1)\psi(n) + \log(n) - \log\Gamma(n) \\
                                      & \leq \frac{1}{2} \log\bigg(\frac{n\mathcal{V}^2}{(\sigma_j^0)^2}\bigg) + \frac{(\sigma_j^0)^2}{2n\mathcal{V}^2} + \frac{(\mu_j^0)^2}{2\mathcal{V}^2} - \frac{1}{2} + \log\bigg(\frac{(\sigma_j^0)^2}{\gamma^2}\bigg) + \frac{\gamma^2}{(\sigma_j^0)^2}-n\\
                & \hspace{2cm}  + n\psi(n) + \log(n) - \log(n-1)! \\   
                & \leq \frac{1}{2} \log\bigg(\frac{n\mathcal{V}^2}{(\sigma_j^0)^2}\bigg) + \frac{(\sigma_j^0)^2}{2n\mathcal{V}^2} + \frac{(\mu_j^0)^2}{2\mathcal{V}^2} - \frac{1}{2} + \log\bigg(\frac{(\sigma_j^0)^2}{\gamma^2}\bigg) + \frac{\gamma^2}{(\sigma_j^0)^2}-n\\
                & \hspace{2cm}  + \bigg( n\log(n)-\frac{n}{2n}-\frac{n}{12n^2}+\frac{n}{120n^4} \bigg) + \log(n) \\
                & \hspace{2cm} +\bigg( - \frac{1}{2}\log(2\pi) +n-1 -n\log(n-1) +\frac{1}{2} \log(n-1) \bigg) \\
                & = \frac{1}{2} \log\bigg(\frac{n\mathcal{V}^2}{(\sigma_j^0)^2}\bigg) + \frac{(\sigma_j^0)^2}{2n\mathcal{V}^2} + \frac{(\mu_j^0)^2}{2\mathcal{V}^2} - \frac{3}{2} + \log\bigg(\frac{(\sigma_j^0)^2}{\gamma^2}\bigg) + \frac{\gamma^2}{(\sigma_j^0)^2}\\
                & \hspace{2cm}  + n\log(\frac{n}{n-1})-\frac{1}{2} + \log(n) - \frac{1}{2}\log(2\pi) +\frac{1}{2} \log(n-1) \\
                & \leq \frac{(\sigma_j^0)^2}{2n\mathcal{V}^2} + \frac{(\mu_j^0)^2}{2\mathcal{V}^2} - \frac{3}{2} + \bigg( \frac{1}{2} \log\bigg(\frac{(\sigma_j^0)^2}{\gamma^4}\bigg) + \frac{\gamma^2}{(\sigma_j^0)^2}  - \frac{1}{2}\log(2\pi) \bigg)\\
                & \hspace{2cm}  + 2 -\frac{1}{2} + \bigg( \frac{1}{2} \log\big(n\mathcal{V}^2\big) + \frac{3}{2} \log(n) \bigg)\\
                & = \frac{(\sigma_j^0)^2}{2n\mathcal{V}^2} + \frac{(\mu_j^0)^2}{2\mathcal{V}^2} + \bigg( \frac{1}{2} \log\bigg(\frac{(\sigma_j^0)^2}{\gamma^4}\bigg) + \frac{\gamma^2}{(\sigma_j^0)^2}  - \frac{1}{2}\log(2\pi) \bigg) + 2 \log(n\sqrt{\mathcal{V}})\\
                & = n \times \frac{1}{n} \bigg[ 2 \log(n\sqrt{\mathcal{V}}) + \frac{(\sigma_j^0)^2}{2n\mathcal{V}^2} + \frac{(\mu_j^0)^2}{2\mathcal{V}^2} + \frac{1}{2} \log\bigg(\frac{(\sigma_j^0)^2}{\gamma^4}\bigg) + \frac{\gamma^2}{(\sigma_j^0)^2}  - \frac{1}{2}\log(2\pi) \bigg] \\
                                     & \leq n R_{j,n}
\end{aligned}$

\vspace{0.2cm}
with $R_{j,n}=\frac{1}{n} \bigvee \frac{1}{n} \bigg[ 2 \log(n\sqrt{\mathcal{V}}) + \frac{(\sigma_j^0)^2}{2n\mathcal{V}^2} + \frac{(\mu_j^0)^2}{2\mathcal{V}^2} + \frac{1}{2} \log\bigg(\frac{(\sigma_j^0)^2}{\gamma^4}\bigg) + \frac{\gamma^2}{(\sigma_j^0)^2}  - \frac{1}{2}\log(2\pi) \bigg].$

\vspace{0.2cm}

The end of the proof is the same as the one used in the Normal-Inverse-Gamma case.

\end{proof}

\vspace{0.2cm}
\subsection{Proof of Theorem \ref{thm-misspecified}}
\vspace{0.2cm}

\begin{proof}
We assume that $\Theta_K(r_{n,K})$ is not empty (otherwise, this is obvious). Applying Theorem 2.7 in \cite{Tempered} for any $\alpha \in (0,1)$, $\theta^* \in \Theta_K(r_{n,K})$:

\begin{equation*}
\begin{split}
\mathbb{E} \bigg[ \int D_{\alpha}( P_{\theta},& P^0 ) \tilde{\pi}_{n,\alpha}(d\theta|X_1^n) \bigg] \leq \frac{\alpha}{1-\alpha} \mathcal{K}(P^0,P_{\theta^*}) \\
 &  + \inf_{\rho \in \mathcal{F}} \bigg\{ \frac{\alpha}{1-\alpha} \int \mathbb{E} \bigg[ \log \frac{P_{\theta^*}(X)}{P_{\theta}(X)} \bigg] \rho(d\theta) + \frac{\mathcal{K}(\rho_p,\pi_p) + \sum_{j=1}^K \mathcal{K}(\rho_j,\pi_j)}{n(1-\alpha)} \bigg\}.
\end{split}
\end{equation*}

\vspace{0.2cm}
Let us take $\rho_{j,n}$ and $\mathcal{A}_{n,K}$ from the definition of $\Theta_K(r_{n,K})$, and $\rho_{p,n}(dp) \propto \mathbf{1}(p \in \mathcal{A}_{n,K}) \pi_p(dp)$:

\begin{equation*}
\begin{split}
\mathbb{E} \bigg[ \int D_{\alpha}( P_{\theta},& P^0 ) \tilde{\pi}_{n,\alpha}(d\theta|X_1^n) \bigg] \leq \frac{\alpha}{1-\alpha} \mathcal{K}(P^0,P_{\theta^*}) \\
 &  + \frac{\alpha}{1-\alpha} \int \mathbb{E} \bigg[ \log \frac{P_{\theta^*}(X)}{P_{\theta}(X)} \bigg] \rho_{p,n}(dp) \prod_{j=1}^K \rho_{j,n}(d\theta_j) + \frac{\mathcal{K}(\rho_{p,n},\pi_p) + \sum_{j=1}^K \mathcal{K}(\rho_{j,n},\pi_j)}{n(1-\alpha)}.
\end{split}
\end{equation*}

\vspace{0.2cm}
We have $\mathcal{K}(\rho_{p,n},\pi_p)=-\log(\pi_p(\mathcal{A}_{n,K}))\leq nKr_{n,K}$ and $\mathcal{K}(\rho_{j,n},\pi_j) \leq nr_{n,K}$ for each $j$ by definition of $\Theta_K(r_{n,K})$. Moreover, using the same argument contained in the proof of Lemma \ref{Do}:

\vspace{0.2cm}
$\begin{aligned}[t]
   \log \frac{P_{\theta^*}(X)}{P_{\theta}(X)} & = \frac{1}{P_{\theta^*}(X)} P_{\theta^*}(X) \log \frac{P_{\theta^*}(X)}{P_{\theta}(X)} \\
                                     & \leq \frac{1}{P_{\theta^*}(X)} \sum\limits_{j=1}^K p_j^* q_{\theta_j^*}(X) \log \frac{p_j^* q_{\theta_j^*}(X)}{p_j q_{\theta_j}(X)} \\
                                     & = \sum\limits_{j=1}^K \frac{p_j^* q_{\theta_j^*}(X)}{P_{\theta^*}(X)} \log \frac{p_j^* }{p_j} + \sum\limits_{j=1}^K \frac{p_j^* q_{\theta_j^*}(X)}{P_{\theta^*}(X)} \log \frac{q_{\theta_j^*}(X)}{q_{\theta_j}(X)} \\
                                     & \leq \sum\limits_{j=1}^K \frac{p_j^* q_{\theta_j^*}(X)}{P_{\theta^*}(X)} \log \frac{p_j^* }{p_j} + \sum\limits_{j=1}^K \log \frac{q_{\theta_j^*}(X)}{q_{\theta_j}(X)}
\end{aligned}$

\noindent and thus, as the support of $\rho_{p,n}$ is on $\mathcal{A}_{n,K}$ where $\log \frac{p_j^* }{p_j} \leq Kr_{n,K}$,

\vspace{0.2cm}
$\begin{aligned}[t]
   \int \mathbb{E} \bigg[ \log \frac{P_{\theta^*}(X)}{P_{\theta}(X)} \bigg] \rho_{p,n}(dp) \prod_{j=1}^K \rho_{j,n}(d\theta_j) & \leq \int \mathbb{E} \bigg[  \sum\limits_{j=1}^K \frac{p_j^* q_{\theta_j^*}}{P_{\theta^*}} \bigg]  \log \frac{p_j^* }{p_j} \rho_{p,n}(dp) \\
                                     &  \hspace{0.5cm} + \sum\limits_{j=1}^K \int \mathbb{E} \bigg[ \log \frac{q_{\theta_j^*}(X)}{q_{\theta_j}(X)} \bigg] \rho_{j,n}(d\theta_j) \\
                                     & \leq \int \mathbb{E} \bigg[ \sum\limits_{j=1}^K \frac{p_j^* q_{\theta_j^*}}{P_{\theta^*}} \bigg] Kr_{n,K} \rho_{p,n}(dp)+ K r_{n,K} \\
                                     & = 2Kr_{n,K}
\end{aligned}$

\vspace{0.2cm}
\noindent which ends the proof as it holds for any $\theta^* \in \Theta_K(r_{n,K})$.

\end{proof}

\vspace{0.2cm}
\subsection{Proof of Corollary \ref{thm-misspecified-Gaussian}}
\vspace{0.2cm}

\begin{proof}
It is sufficient to show that $\mathcal{S}_K \times [-L,L]^K \subset \Theta_K(r_{n,K})$ for 
$$
r_{n,K}=\frac{4\log(nK)}{n}  \bigvee_{j=1}^{K} \frac{1}{n} \bigg[ \frac{1}{2} \log\bigg(\frac{n}{2}\bigg) + \frac{1}{n\mathcal{V}^2} + \log\big({\mathcal{V}}\big) + \frac{L^2}{2\mathcal{V}^2} - \frac{1}{2} \bigg],
$$
the stated oracle inequality is a direct corollary of Theorem \ref{thm-misspecified}. For that, let us take any $\theta^* \in \mathcal{S}_K \times [-L,L]^K$ and show that it satisfies the conditions in the definition of $\Theta_K(r_{n,K})$.

\vspace{0.2cm}
The existence of a set $\mathcal{A}_{n,K}$ fulfilling the first condition has already been done in the proof of Lemma \ref{lemma-Dirichlet} as $\frac{4\log(nK)}{n} \leq r_{n,K}$.

\vspace{0.2cm}
We define distributions $\rho_{j,n} \in \mathcal{M}_1^+(\Theta)$ by Gaussians of mean $\theta_j^*$ and variance $\frac{2}{n}$ ($j=1,...,K$) and we show that for $j=1,...,K$:
\begin{equation*}
  \int \mathbb{E}\bigg[\log\bigg(\frac{q_{\theta_{j}^*}(X)}{q_{\theta_{j}}(X)} \bigg)\bigg]  \rho_{j,n}(d\theta_{j}) \leq r_{n,K} \hspace{0.2cm},  \hspace{0.5cm}
  \mathcal{K}(\rho_{j,n},\pi_{j}) \leq n r_{n,K}.
\end{equation*}
We start from
\[
\log\bigg(\frac{q_{\theta_j^*}(X)}{q_{\theta_j}(X)}\bigg) = \frac{(\theta_j-\theta_j^*)^2}{2} - {(X-\theta_j^*)(\theta_j-\theta_j^*)}
\]
and if we take the mean of this quantity with respect to $P^0$, we obtain:
\[
\mathbb{E}\bigg[\log\bigg(\frac{q_{\theta_{j}^*}(X)}{q_{\theta_{j}}(X)} \bigg)\bigg] = \frac{(\theta_j-\theta_j^*)^2}{2} - {(\mathbb{E}X-\theta_j^*)(\theta_j-\theta_j^*)}
\]
and as $\theta_j-\theta_j^*$ is a zero-mean random variable, we have:

\vspace{0.2cm}
$\begin{aligned}[t]
   \int \mathbb{E}\bigg[\log\bigg(\frac{q_{\theta_{j}^*}(X)}{q_{\theta_{j}}(X)} \bigg)\bigg] \rho_{j,n}(d\theta_j) & = \frac{1}{2}  \mathbb{E}_{\theta_j\sim\rho_{j,n}}[(\theta_j-\theta_j^*)^2] - {(\mathbb{E}X-\theta_j^*)} \mathbb{E}_{\theta_j\sim\rho_{j,n}}[\theta_j-\theta_j^*]  \\
                                     & = \frac{1}{2} \times \frac{2}{n} \\
                                     & \leq r_{n,K}.
\end{aligned}$

\vspace{0.2cm}
Then, we conclude according to Lemma \ref{thm-kl-divergence-gaussians}: 

\vspace{0.2cm}
$\begin{aligned}[t]
   \mathcal{K}(\rho_{j,n},\pi_j) & = \frac{1}{2} \log\bigg(\frac{n\mathcal{V}^2}{2}\bigg) + \frac{1}{n\mathcal{V}^2} + \frac{(\theta_j^*)^2}{2\mathcal{V}^2} - \frac{1}{2}  \\
                                     & = \frac{1}{2} \log\bigg(\frac{n}{2}\bigg) + \frac{1}{n\mathcal{V}^2} + \log\big({\mathcal{V}}\big) + \frac{(\theta_j^*)^2}{2\mathcal{V}^2} - \frac{1}{2} \\
                                     & \leq n \times \frac{1}{n} \bigg[ \frac{1}{2} \log\bigg(\frac{n}{2}\bigg) + \frac{1}{n\mathcal{V}^2} + \log\big({\mathcal{V}}\big) + \frac{L^2}{2\mathcal{V}^2} - \frac{1}{2} \bigg] \\
                                     & \leq n r_{n,K}.
\end{aligned}$

\end{proof} 

\vspace{0.2cm}
\subsection{Proof of Theorem \ref{thm-model-selection}}
\vspace{0.2cm}

Here, we cannot directly use the results from \cite{Tempered}. So we prove this theorem from scratch, by following the main steps outlined in \cite{bhattacharya2016bayesian,Tempered}  with some adaptation.

\begin{proof}
For any $\alpha \in (0,1)$ and $\theta \in \Omega$, by definition of the Renyi divergence and using $D_\alpha(P^{\otimes n},R^{\otimes n})=nD_\alpha(P,R)$ as data are i.i.d.:
$$
\mathbb{E}\bigg[ \exp\bigg(-\alpha r_{n}(P_\theta,P^0) + (1-\alpha)n D_\alpha(P_\theta,P^0)\bigg) \bigg] = 1
$$
Thus, integrating and using Fubini's theorem, $$ 
\mathbb{E}\bigg[ \int \exp\bigg(-\alpha r_{n}(P_\theta,P^0) + (1-\alpha)n D_\alpha(P_\theta,P^0) \bigg) \pi(d\theta) \bigg] = 1
$$
Using Lemma \ref{thm-dv},
$$ 
\mathbb{E}\bigg[ \exp\bigg( \sup_{\rho \in \mathcal{M}_1^+(\Omega)} \bigg\{ \int \bigg( -\alpha r_{n}(P_\theta,P^0) + (1-\alpha)n D_\alpha(P_\theta,P^0) \bigg) \rho(d\theta) - \mathcal{K}(\rho,\pi) \bigg\} \bigg) \bigg] = 1.
$$
Note that \cite{bhattacharya2016bayesian,Tempered} also used Lemma \ref{thm-dv} in their proofs, this is inspired by the PAC-Bayesian theory \cite{catoni2004statistical,MR2483528}. It is interesting to note that Lemma \ref{thm-dv} is at the core of VB: it is used to provide approximation algorithms, and also to prove the consistency of VB.
Thanks to Jensen's inequality, $$ 
\mathbb{E}\bigg[ \sup_{\rho \in \mathcal{M}_1^+(\Omega)} \bigg\{ \int \bigg( -\alpha r_{n}(P_\theta,P^0) + (1-\alpha)n D_\alpha(P_\theta,P^0) \bigg) \rho(d\theta) - \mathcal{K}(\rho,\pi) \bigg\} \bigg] \leq 0
$$
Therefore, when considering $\tilde{\pi}^{\hat{K}}_{n,\alpha}(.|X_1^n)$ as a distribution on $\mathcal{M}_1^+(\Omega)$ with all its mass on $\Theta_{\hat{K}}$, $$ 
\mathbb{E}\bigg[ \int \bigg( -\alpha r_{n}(P_\theta,P^0) + (1-\alpha)n D_\alpha(P_\theta,P^0) \bigg) \tilde{\pi}_{n,\alpha}^{\hat{K}}(d\theta|X_1^n) - \mathcal{K}(\tilde{\pi}^{\hat{K}}_{n,\alpha}(.|X_1^n),\pi) \bigg] \leq 0
$$
Using $\mathcal{K}(\tilde{\pi}^{\hat{K}}_{n,\alpha}(.|X_1^n),\pi)=\mathcal{K}(\tilde{\pi}^{\hat{K}}_{n,\alpha}(.|X_1^n),\Pi_{\hat{K}})+\log(\frac{1}{\pi_{\hat{K}}})$, we rearrange terms:
\begin{multline*}
\mathbb{E}\bigg[ \int D_\alpha(P_\theta,P^0) \tilde{\pi}_{n,\alpha}^{\hat{K}}(d\theta|X_1^n) \bigg] \\ \leq 
\mathbb{E}\bigg[ \frac{\alpha}{1-\alpha} \int \frac{r_{n}(P_\theta,P^0)}{n} \tilde{\pi}_{n,\alpha}^{\hat{K}}(d\theta|X_1^n) + \frac{\mathcal{K}(\tilde{\pi}^{\hat{K}}_{n,\alpha}(.|X_1^n),\Pi_{\hat{K}})}{n(1-\alpha)}+\frac{\log(\frac{1}{\pi_{\hat{K}}})}{n(1-\alpha)} \bigg] 
\end{multline*}
Thus, by definition of $\hat{K}$,
\begin{multline*}
\mathbb{E}\bigg[ \int D_\alpha(P_\theta,P^0) \tilde{\pi}_{n,\alpha}^{\hat{K}}(d\theta|X_1^n) \bigg] \\ \leq \mathbb{E}\bigg[  \inf_{K \geq 1} \bigg\{ \frac{\alpha}{1-\alpha} \int \frac{r_{n}(P_\theta,P^0)}{n} \tilde{\pi}_{n,\alpha}^{{K}}(d\theta|X_1^n) + \frac{\mathcal{K}(\tilde{\pi}^{{K}}_{n,\alpha}(.|X_1^n),\Pi_{{K}})}{n(1-\alpha)}+\frac{\log(\frac{1}{\pi_{{K}}})}{n(1-\alpha)} \bigg\} \bigg]
\end{multline*}
which leads to 
\begin{multline*}
\mathbb{E}\bigg[ \int D_\alpha(P_\theta,P^0) \tilde{\pi}_{n,\alpha}^{\hat{K}}(d\theta|X_1^n) \bigg] \\
\leq \inf_{K \geq 1} \bigg\{ \mathbb{E}\bigg[ \frac{\alpha}{1-\alpha} \int \frac{r_{n}(P_\theta,P^0)}{n} \tilde{\pi}_{n,\alpha}^{{K}}(d\theta|X_1^n) + \frac{\mathcal{K}(\tilde{\pi}^{{K}}_{n,\alpha}(.|X_1^n),\Pi_{{K}})}{n(1-\alpha)}+\frac{\log(\frac{1}{\pi_{{K}}})}{n(1-\alpha)} \bigg] \bigg\}
\end{multline*}
and by definition of $\tilde{\pi}^{{K}}_{n,\alpha}(.|X_1^n)$, \begin{multline*}
\mathbb{E}\bigg[ \int D_\alpha(P_\theta,P^0) \tilde{\pi}_{n,\alpha}^{\hat{K}}(d\theta|X_1^n) \bigg] \\ \leq \inf_{K \geq 1} \bigg\{ \mathbb{E}\bigg[ \inf_{\rho \in \mathcal{M}_1^+(\Theta_K)} \bigg\{ \frac{\alpha}{1-\alpha} \int \frac{r_{n}(P_\theta,P^0)}{n} \rho(d\theta) + \frac{\mathcal{K}(\rho,\Pi_{{K}})}{n(1-\alpha)} \bigg\} +\frac{\log(\frac{1}{\pi_{{K}}})}{n(1-\alpha)} \bigg] \bigg\}.
\end{multline*}
Then, \begin{multline*}
\mathbb{E}\bigg[ \int D_\alpha(P_\theta,P^0) \tilde{\pi}_{n,\alpha}^{\hat{K}}(d\theta|X_1^n) \bigg]\\  \leq \inf_{K \geq 1} \inf_{\rho \in \mathcal{M}_1^+(\Theta_K)} \bigg\{ \mathbb{E}\bigg[ \frac{\alpha}{1-\alpha} \int \frac{r_{n}(P_\theta,P^0)}{n} \rho(d\theta) + \frac{\mathcal{K}(\rho,\Pi_{{K}})}{n(1-\alpha)} +\frac{\log(\frac{1}{\pi_{{K}}})}{n(1-\alpha)} \bigg] \bigg\}.
\end{multline*}
And finally, $$
\mathbb{E}\bigg[ \int D_\alpha(P_\theta,P^0) \tilde{\pi}_{n,\alpha}^{\hat{K}}(d\theta|X_1^n) \bigg] \leq \inf_{K \geq 1} \inf_{\rho \in \mathcal{M}_1^+(\Theta_K)} \bigg\{ \frac{\alpha}{1-\alpha} \int \mathcal{K}(P^0,P_\theta) \rho(d\theta) + \frac{\mathcal{K}(\rho,\Pi_{{K}})}{n(1-\alpha)} +\frac{\log(\frac{1}{\pi_{{K}}})}{n(1-\alpha)} \bigg\}.
$$
To conclude, we just need to upper bound the function inside the infimum over all integers $K$'s by $\frac{\alpha}{1-\alpha} \inf_{\theta^* \in \Theta_K(r_{n,K})} \mathcal{K}(P^0,P_{\theta^*}) + \frac{1+\alpha}{1-\alpha}2Kr_{n,K} +  \frac{\log(\frac{1}{\pi_K})}{n(1-\alpha)}$. This is direct: if the set $\Theta_K(r_{n,K})$ is not empty (otherwise the inequality is obvious) we notice that $\mathcal{K}(P^0,P_{\theta}) = \mathcal{K}(P^0,P_{\theta^*}) + \mathbb{E} \bigg[ \log \frac{P_{\theta^*}(X_i)}{P_{\theta}(X_i)} \bigg]$ for any $\theta^* \in \Theta_K(r_{n,K})$ and then we follow the sketch of the proof of Theorem \ref{thm-misspecified}.

\end{proof}

\vspace{0.2cm}
\subsection{Algorithms}
\vspace{0.2cm}

We now provide the derivations leading to the algorithms described in the paper.

\vspace{0.2cm}
\subsubsection{Algorithm 1}
\vspace{0.2cm}

We apply a coordinate descent on variables $\omega^1 \in \mathcal{S}_K$,..., $\omega^n \in \mathcal{S}_K$, $\rho_p \in \mathcal{M}_1^+(\mathcal{S}_K)$, $\rho_1 \in \mathcal{M}_1^+(\Theta)$,..., and $\rho_K \in \mathcal{M}_1^+(\Theta)$ in order to solve the optimization program:
\begin{equation*}
\begin{split}
\min_{\rho \in \mathcal{F}, \hspace{0.1cm} w \in \mathcal{S}_K^n} \bigg\{ - \alpha \sum\limits_{i=1}^n \sum\limits_{j=1}^K \omega_j^i & \bigg( \int \log(p_j) \rho_p(dp) + \int \log(q_{\theta_j}(X_i)) \rho_j(d\theta_j) \bigg) \\
& + \alpha \sum\limits_{i=1}^n \sum\limits_{j=1}^K \omega_j^i \log(\omega_j^i)
+ \mathcal{K}(\rho_p,\pi_p) + \sum\limits_{j=1}^K \mathcal{K}(\rho_j,\pi_j) \bigg\}.
\end{split} 
\end{equation*}

\vspace{0.2cm}
We explain how to obtain Algorithm 1.

\vspace{0.2cm}
\subsubsection*{Optimization with respect to $\omega^i \in \mathcal{S}_K$:}
\vspace{0.2cm}

First, we fix $\omega^{\ell} \in \mathcal{S}_K$ for $\ell \ne i$, $\rho_p \in \mathcal{M}_1^+(\mathcal{S}_K)$ and $\rho_j \in \mathcal{M}_1^+(\Theta)$ for $j=1,...,K$, and we solve the program with respect to $\omega^i \in \mathcal{S}_K$, which becomes:
\begin{equation*}
\begin{split}
\min_{\omega^i \in \mathcal{S}_K} \bigg\{ \sum\limits_{j=1}^K \omega_j^i & \bigg( \log(\omega_j^i) - \int \log(p_j) \rho_p(dp) - \int \log(q_{\theta_j}(X_i)) \rho_j(d\theta_j) \bigg) \bigg\}.
\end{split} 
\end{equation*}
Put $\textbf{E}=\{1,...,K\}$, $\lambda=\left(\frac{1}{K},...,\frac{1}{K}\right)$ and $h(j)=\int \log(p_j) \rho_p(dp) + \int \log(q_{\theta_j}(X_i)) \rho_j(d\theta_j)$ and use Lemma \ref{thm-dv} to obtain:
$$
w^i_j \propto \exp\bigg( \int \log(p_j) \rho_p(dp) + \int \log(q_{\theta_j}(X_i)) \rho_j(d\theta_j) \bigg).
$$

\vspace{0.2cm}
\subsubsection*{Optimization with respect to $\rho_p \in \mathcal{M}_1^+(\mathcal{S}_K)$:}
\vspace{0.2cm}
Now, we fix $\omega^{i} \in \mathcal{S}_K$ for $i=1,...,n$, and $\rho_j \in \mathcal{M}_1^+(\Theta)$ for $j=1,...,K$, and we solve the program with respect to $\rho_p \in \mathcal{M}_1^+(\mathcal{S}_K)$, which becomes:
\begin{equation*}
\begin{split}
\min_{\rho_p \in \mathcal{M}_1^+(\mathcal{S}_K)} \bigg\{ - \alpha \sum\limits_{i=1}^n \sum\limits_{j=1}^K \omega_j^i \int \log(p_j) \rho_p(dp)
+ \mathcal{K}(\rho_p,\pi_p) \bigg\}.
\end{split} 
\end{equation*}

\noindent Using Lemma \ref{thm-dv} for $\textbf{E}=\mathcal{S}_K$, $\lambda=\pi_p$ and $h(p)=\alpha \sum\limits_{i=1}^n \sum\limits_{j=1}^K \omega_j^i \log(p_j)$, we get directly the solution: 
$$
\rho_p(dp) \propto \exp\bigg( \alpha \sum\limits_{i=1}^n \sum\limits_{j=1}^K \omega^i_j \log(p_j) \bigg) \pi_p(dp).
$$

\vspace{0.2cm}
\subsubsection*{Optimization with respect to $\rho_j \in \mathcal{M}_1^+(\Theta)$:}
\vspace{0.2cm}
Now, we fix $\omega^{i} \in \mathcal{S}_K$ for $i=1,...,n$, $\rho_p \in \mathcal{M}_1^+(\mathcal{S}_K)$ and $\rho_\ell \in \mathcal{M}_1^+(\Theta)$ for $\ell \ne j$, and we solve the program with respect to $\rho_j \in \mathcal{M}_1^+(\Theta)$, which becomes:
\begin{equation*}
\begin{split}
\min_{\rho_j \in \mathcal{M}_1^+(\Theta)} \bigg\{ - \alpha \sum\limits_{i=1}^n \omega_j^i \int \log(q_{\theta_j}(X_i)) \rho_j(d\theta_j)
+ \mathcal{K}(\rho_j,\pi_j) \bigg\}.
\end{split} 
\end{equation*}

\noindent Using Lemma \ref{thm-dv} for $\textbf{E}=\Theta$, $\lambda=\pi_j$ and $h(\theta_j)=\alpha \sum\limits_{i=1}^n \omega_j^i \log(q_{\theta_j}(X_i))$, we get directly the solution: 
$$
\rho_j(d\theta_j) \propto \exp\bigg( \alpha \sum\limits_{i=1}^n \omega^i_j \log(q_{\theta_j}(X_i)) \bigg) \pi_j(d\theta_j)
$$

\vspace{0.2cm}
\subsubsection{Application to multinomial mixture models}
\vspace{0.2cm}

We simply use 
$$
\int \log(p_j) \rho_p(dp) = \mathbb{E}_{p\sim \rho_p}[\log(p_j)] =  \psi(\phi_j)-\psi(\sum\limits_{\ell=1}^K \phi_{\ell}),
$$
$$
\int \log(q_{\theta_j}(X_i)) \rho_j(d\theta_j) =  \mathbb{E}_{\theta_j \sim \rho_j}[\log(\theta_{X_i,j})] = \psi(\gamma_{X_i,j}) - \psi\big(\sum\limits_{v=1}^V \gamma_{vj}\big),
$$
$$
\exp\bigg( \alpha \sum\limits_{i=1}^n \sum\limits_{j=1}^K \omega^i_j \log(p_j) \bigg) \pi_p(p) \propto \prod_{j=1}^K p_j^{\alpha_j+\alpha \sum_{i=1}^n \omega^i_j-1},
$$
$$
\exp\bigg( \alpha \sum\limits_{i=1}^n \omega^i_j \log(q_{\theta_j}(X_i)) \bigg) \pi_j(\theta_j) \propto \prod_{v=1}^V \theta_{vj}^{\beta_v+\alpha \sum\limits_{i=1}^n \omega_{j}^i \mathds{1}(X_i=v)-1}.
$$
We recognize a Dirichlet distribution.

\vspace{0.2cm}
\subsubsection{Application to Gaussian mixture models}
\vspace{0.2cm}

For Gaussian mixtures, use
$$
\int \log(p_j) \rho_p(dp) = \mathbb{E}_{p\sim \rho_p}[\log(p_j)] =  \psi(\phi_j)-\psi(\sum\limits_{\ell=1}^K \phi_{\ell}),
$$
$$
\int \log(q_{\theta_j}(X_i)) \rho_j(d\theta_j) = - \frac{1}{2} \mathbb{E}_{\theta_j \sim \rho_j}[(\theta_{j}-X_i)^2] + \textnormal{cst} = - \frac{1}{2}\big\{ s_j^2 + (n_j-X_i)^2 \big\} + \textnormal{cst},
$$
$$
\exp\bigg( \alpha \sum\limits_{i=1}^n \sum\limits_{j=1}^K \omega^i_j \log(p_j) \bigg) \pi_p(p) \propto \prod_{j=1}^K p_j^{\alpha_j+\alpha \sum_{i=1}^n \omega^i_j-1},
$$
\begin{equation*}
\begin{split}
\exp\bigg( \alpha \sum\limits_{i=1}^n \omega^i_j \log(q_{\theta_j}(X_i)) \bigg) \pi_j(\theta_j) & \propto \exp\bigg( - \frac{\alpha}{2} \sum\limits_{i=1}^n \omega^i_j (\theta_{j}-X_i)^2 \bigg) \exp\bigg( - \frac{1}{2\mathcal{V}^2} \theta_j^2 \bigg) \\
& \propto \exp\bigg( - \frac{{1/\mathcal{V}^2+\alpha \sum_{i=1}^n \omega_{j}^i}}{2} \bigg(\theta_{j}-\frac{\alpha \sum_{i=1}^n \omega_{j}^iX_i}{1/\mathcal{V}^2+\alpha \sum_{i=1}^n \omega_{j}^i}\bigg)^2 \bigg) .
\end{split} 
\end{equation*}
We recognize a Gaussian distribution.

\newpage

\subsection*{Supplementary material}

\vspace{0.4cm}
We provide in this supplementary material a very short simulation study. Our objective is not to compare extensively EM to CAVI as this was already done in many papers (mentioned in the main body of the paper). We just show on a low-dimensional example that the properties of VB with $\alpha=1/2$ and $\alpha=1$ (CAVI) are very similar to each other, and also to EM.

\vspace{0.2cm}
We compare our algorithm for $\alpha=0.5$ and $\alpha=1$ (equivalent to CAVI) to EM algorithm for unit-variance Gaussian mixture parameters estimation. We consider 10 different unit-variance Gaussian mixtures, where the parameters $(p^0,\theta_1^0,\theta_2^0,\theta_3^0)$ are generated independently from a Dirichlet distribution $p^0 \sim \mathcal{D}_K(2/3,2/3,2/3)$ and Gaussians $\theta_j^0 \sim \mathcal{N}(0,10)$ for $j=1,2,3$. From these mixtures, we create 10 different datasets which contain 1000 i.i.d. realizations of the corresponding mixtures. We compare our algorithms using the Mean Average Error (MAE) between the estimates and the true parameters. For each dataset, we run each algorithm 5 times and keep the one with the lowest MAE in order to avoid situations where the initialization leads to a local optimum. Then, we average the resulting MAEs over the different datasets to obtain the final values of the MAE. We also record the standard deviation of the MAE over the different datasets. The following table summarizes the results. Values in parenthesis represent the standard deviations of the computed MAEs, and the three components are ordered in ascending values. The three procedures are comparable both in terms of estimation precision and computational efficiency :

\vspace{0.4cm}

\begin{center}
\centering

\begin{tabular}{|c|c|c|c|c|c|c|c|c|c|c|}

\hline
Algorithm & $p$ & $\theta_1$ & $\theta_2$ & $\theta_3$ \\
\hline
VB ($\alpha=0.5$) & 0.033 (0.020) & 0.137 (0.297) & 0.383 (1.108) & 0.054 (0.047) \\
\hline
VB ($\alpha=1$) & 0.033 (0.020) & 0.139 (0.207) & 0.364 (0.968) & 0.056 (0.039) \\
\hline
EM & 0.033 (0.021) & 0.141 (0.219) & 0.364 (0.968) & 0.059 (0.047) \\
\hline

\end{tabular}
\end{center}

The notebook is available on the second author webpage:
\begin{verbatim}
 http://alquier.ensae.net/packages.html
\end{verbatim}

\end{document}